\documentclass[11pt]{amsart}
\usepackage{fullpage}
\usepackage{geometry}                % See geometry.pdf to learn the layout options. There are lots.
\usepackage{graphicx}
\usepackage{amssymb}
\usepackage{epstopdf}
\usepackage{color}
\usepackage{bbm, dsfont}
\newtheorem{theorem}{Theorem}
\newtheorem{lemma}{Lemma}[section]

\newtheorem{proposition}{Proposition}[section]

\theoremstyle{definition}
\newtheorem{definition}{Definition}
\theoremstyle{remark}
\newtheorem{remark}{Remark}

\newtheorem{notation}{Notation}
\theoremstyle{example}

\usepackage{latexsym,amssymb,amsmath,amsfonts,amsthm}
\usepackage{amsthm,amsmath}
\usepackage{float}
\usepackage{enumerate}
\usepackage{epic}
\usepackage{fullpage}
\usepackage{bbm}
\usepackage{subfigure}
\usepackage{graphicx}
\usepackage[toc,page]{appendix}
\usepackage{amsfonts}
\usepackage[all]{xy}
\usepackage{caption}
 \usepackage{geometry}              
\usepackage{graphicx}
\usepackage{enumerate}
\usepackage{amssymb}
\usepackage{epstopdf}
\usepackage{color}
\usepackage{bbm, dsfont}
\usepackage{hyperref}
\usepackage[utf8]{inputenc}
\usepackage{amssymb,amsthm,amsmath,amssymb,wrapfig,dsfont}
\usepackage[dvipsnames]{xcolor}
\definecolor{myred}{RGB}{251,154,133}
\definecolor{myblue}{RGB}{153,206,227}
\definecolor{mylightblue}{RGB}{0, 150, 255}
\definecolor{mygreen}{RGB}{32, 210, 64}
\definecolor{mygray}{RGB}{220, 220, 220}

\usepackage{tikz}
\usetikzlibrary{decorations.pathmorphing}
\tikzset{snake it/.style={decorate, decoration=snake}}
\usetikzlibrary{shapes.geometric,positioning,decorations.pathreplacing}

\DeclareFontFamily{OML}{rsfs}{\skewchar\font'177}
\DeclareFontShape{OML}{rsfs}{m}{n}{ <5> <6> rsfs5 <7> <8> <9>
rsfs7 <10> <10.95> <12> <14.4> <17.28> <20.74> <24.88> rsfs10 }{}
\DeclareMathAlphabet{\mathfs}{OML}{rsfs}{m}{n}

\newcommand{\BE}{{\mathbb{E}}}

\newcommand{\BH}{{\mathbb{H}}}

\newcommand{\BL}{{\mathbb{L}}}

\newcommand{\BZ}{{\mathbb{Z}}}

\newcommand{\CF}{{\mathcal{F}}}

\newcommand{\CH}{{\mathcal{H}}}

\newcommand{\CP}{{\mathcal{P}}}

\newcommand{\prob}{{\bf P}}
\newcommand{\p}{{\bf P}}
\newcommand{\e}{{\bf E}}
\newcommand{\bae}{\begin{equation}\begin{aligned}}
\newcommand{\eae}{\end{aligned}\end{equation}}

\newcommand{\pin}{\partial^{in}}
\newcommand{\pout}{\partial^{out}}

\newcommand{\ep}{{\epsilon}}

\newcommand{\harm}{\CH}

\newcommand{\sharm}{\CH^s}
\newcommand{\eharm}{\CH^e}

\DeclareFontFamily{OML}{rsfs}{\skewchar\font'177}
\DeclareFontShape{OML}{rsfs}{m}{n}{ <5> <6> rsfs5 <7> <8> <9>
	rsfs7 <10> <10.95> <12> <14.4> <17.28> <20.74> <24.88> rsfs10 }{}
\DeclareMathAlphabet{\mathfs}{OML}{rsfs}{m}{n}

\usepackage{amsmath}
\usepackage{mathrsfs}

\begin{document}
\title{Scaling limit of DLA on a long line segment}

\author{Yingxin Mu}
\address[Yingxin Mu]{Peking University}
\email{muyingxin@pku.edu.cn}

\author{Eviatar B. Procaccia}
\address[Eviatar B. Procaccia]{Texas A\&M University}
\urladdr{www.math.tamu.edu/~procaccia}
\email{eviatarp@gmail.com}

\author{Yuan Zhang}
\address[Yuan Zhang]{Peking University}
\email{zhangyuan@math.pku.edu.cn}

\thanks{The authors would like to thank an anonymous fat cat in the Temple of Great Enlightenment (Dajue Si)}

\maketitle
%\tableofcontents
\begin{abstract}
	In this paper, we prove that the bulk of DLA starting from a long line segment on the $x$-axis has a scaling limit to the stationary DLA process (SDLA). The main phenomenological difficulty is the multi-scale, non-monotone interaction of the DLA arms. We overcome this via a coupling scheme between the two processes and an intermediate DLA process with absorbing mesoscopic boundary segments.  
	
%	the infinite SDLA starting from the $x$-axis can be seen as the proper scaling limit of the edge DLA starting from a proportional line segment , with the help of the intermediate processes we constructed.
\end{abstract}

% !TEX root = Convergences_to_SDLA_3.12.19.tex

\section{Introduction}
\begin{figure}[h!]
\includegraphics[width=5cm]{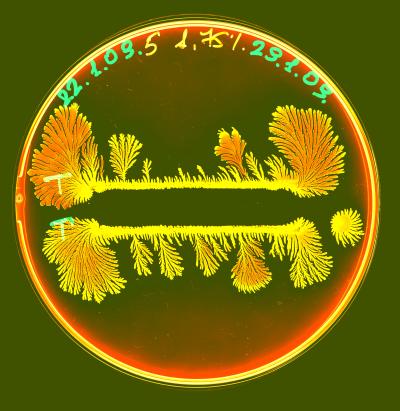}
\caption{Competing bacterial colonies: picture produced in the lab of the late Prof. Eshel Ben-Jacob at Tel-Aviv University.\label{fig:bacteria}}
\end{figure}

In this paper, we establish a scaling limit result for the bulk of DLA on $\BZ^2$ starting from a long line segment. The phenomena of a stationary behavior at the bulk was produced in experimental settings such as in the case of competing bacterial growth on a low nutrient medium (See figure \ref{fig:bacteria} and \cite{be2009deadly}). 

\begin{figure}[h!]
\includegraphics[width=11cm]{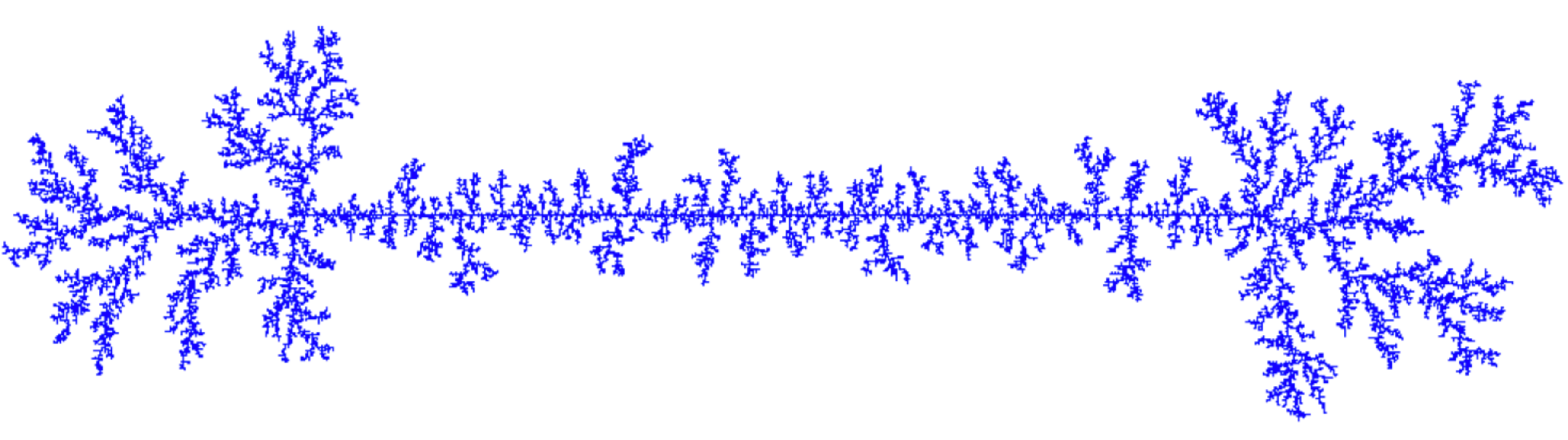}
\caption{A (non-precise) computer simulation of EDLA starting from a long line segment, simulation for qualitative illustration only.\label{fig:edla_line}}
\end{figure}

We consider the edge diffusion limited aggregation (EDLA) on $\BZ^2$, an increasing edge-set process. It grows by adding edges recursively according to the {\bf Edge Harmonic Measure} (the last edge traversed by a random walk coming from infinity before hitting the set). If we start the process from a long line segment, one can observe that in the bulk, the DLA trees tend to grow "upwards" and have similar distribution (See figure \ref{fig:edla_line}). 

In this paper we prove that the bulk of the EDLA starting from a long line segment converges weakly to the infinite stationary DLA (SDLA) process who's existence was established in \cite{procaccia2019stationary}. The SDLA is a continuous time edge-set process on the upper planar lattice generated using a stationary version of the harmonic measure ({\bf stationary harmonic measure}) defined and studied in \cite{procaccia2018stationary,procaccia2017sets,1}.  Several other stationary aggregation processes were recently studied (see \cite{antunovic2017stationary,berger2014stretched}) with some common universal behavior such as a.s. finiteness of all trees.

Before stating the main result, we first need to introduce some terminology.

\subsection{Notations and statement of main results}

Let $\BZ^2$ be the plane square lattice. For any $x=(x(1),x(2))\in \BZ^2$, where $x(1)$ is the first coordinate and $x(2)$ is the second coordinate of $x$, let $\|x\|$ be the $L^2$ norm of vertex $x$. We may turn $\BZ^2$ into a directed graph, by adding a pair of parallel directed edges with opposite orientations between each pair $x,y\in\BZ^2$ with $\|x-y\|=1$. We denote this directed lattice by $\vec\BL^2=(\BZ^2,\vec\BE^2)$ with vertex set $\BZ^2$ and edge set $\vec\BE^2$. For any subset $A\subsetneqq \BZ^2$, intuitively we define $ \vec A$ to be the subgraph of $\vec \BL^2$ whose edge set collects all edges such that both endpoints of these edges are in $A$. Moreover, let $|A|$ be the cardinality of $A$, and if $0\in A$, let
$$
\|A\|=\sup_{x\in A}\|x\|
$$
be the radius of $A$. For any directed edge $ \vec e=x\to y\in \vec\BL^2$, we use $\vec e(1)=x$ and $\vec e(2)=y$ to denote the starting and ending point of $\vec e$.  
 We use
$$
\partial^{in} A=\{x\in A: s.t. \ \exists y\notin A,  \|x-y\|=1\},
$$ 
and
$$
\partial^{out} A=\{x\notin A: s.t. \ \exists y\in A, \|x-y\|=1\}
$$ 
to denote the inner and outer boundaries with respect to vertices. And we use 
$$
\partial^{e} A=\left\{\vec e\in\vec \BL^2: s.t. \ \vec e(1)\in \partial^{out} A, \vec e(2)\in \partial^{in} A\right\}
$$
to denote the edge boundary of $A$ in terms of edges and $\widetilde{\partial^e A}$ to denote the collection of all its inverse edges. Let $\BH$ be the upper half plane. For any $n\ge 0$ we define 
$$
\ell_n=\{(x,n):x\in\BZ\}
$$
as the horizontal line in $\BH$, with $\ell_0$ as the $x-$axis. Moreover, for each $x\in \BZ^2$, let $\prob_x$ be the distribution of the simple random walk $\{S_n\}_{n=0}^\infty$ starting from $x$. And for any $A\subseteq\BZ^2$, one can define the stopping times
$$
\begin{aligned}
&\bar\tau_A=\inf\{n\ge 0: \ S_n\in A\},\\
&\tau_A=\inf\{n\ge 1: \ S_n\in A\}
\end{aligned}
$$
to be the first hitting time and the first returning time respectively. When $A=B(0,R)$, the open ball centered at the origin of radius $R$, we abbreviate them to $\bar\tau_R$ and $\tau_R$.  Here we consider a variant of the DLA model, dubbed edge DLA (EDLA) driven by the 2-dimensional harmonic measure on edges:
\begin{proposition}
\label{prop_harm_out}
For any finite subset $A\subseteq\BZ^2$ and any edge $\vec e$ of $\vec\BL^2$, then the limit
$$
\lim_{\|z\|\to\infty}\prob_{z}\left(\tau_{A}=\tau_{\vec e(2)}, \ S_{\tau_{\vec e(2)}-1}= \vec e(1)\right)
$$
exists. We call the limit above the {\bf Edge Harmonic Measure} of $\vec e$ with respect to $A$, denoted by $\eharm_A(\vec{e})$. 
\end{proposition} 

One may also define the harmonic measure with respect to a vertex $x\in \pout A$ as
$$
\eharm_A(x)=\sum_{\vec e: \ \vec e(1)=x} \eharm_A(\vec e). 
$$
Note that for all $x\in \pin A$, 
$$
\sum_{\vec e: \ \vec e(2)=x} \eharm_A(\vec e)=\harm_A(x)
$$
where $\harm$ stands for the regular harmonic measure on $\BZ^2$. This also implies that 
$$
\sum_{\vec e} \eharm_A(\vec e)=1.
$$
\begin{remark}
However, for $x\in \pout A$, it is important to note that $\eharm_A(x)\not=\harm_A(x)$. 
\end{remark}
With the {\bf Edge Harmonic Measure}, we give a formal description of the EDLA model.
\begin{notation}
Without loss of generality, we often use $V$ and $E$ to distinguish the vertex set from the edge set.
\end{notation}
\begin{definition}\label{def_EDLA}
	\label{def_intermediate_process}
	For any finite $B\subseteq \BZ^2$, one may define the EDLA process $EA^B_t=(EV^B_{t},EE^B_{t})$ to be a continuous time Markov process on the set of all subgraphs of $\vec \BL^2$ such that 
	\begin{itemize}
		\item $EA^B_0=(B,\emptyset)$. 
		\item At any time $t\ge0$, for all edges $\vec e\in \partial^{e} (EV^B_{t-})$, independent Poisson clocks of intensity
		$$
		\lambda(EV^B_{t-},\vec e)=\eharm_{EV^B_{t-}}(\vec e)
		$$
		are placed on $\vec e$. 
		\item If the clock at an edge $\vec e\in \partial^{e} (EV^B_{t-})$ rings at time $t$, let 
		$$
		EA^B_{t}=\left(EV^B_{t-}\cup \{\vec e(1)\},EE^B_{t-}\cup \{\vec e\}\right), 
		$$ 
		and update all the transition rates. 
	\end{itemize}
\end{definition}
\begin{remark}\label{remk1_1}
Note that $EV_t$ forms a vertex-set process which is identically distributed to the Outer DLA process $OA_t$ defined in Definition 1 of \cite{procaccia2019stationary}.
\end{remark}

For any finite $B\subseteq \BZ^2$, the well-definedness of $EA^B_t $ is obvious since the total transition rate is 1.  In this paper, we also use $EA^n_t$ in abbreviation for the case when $EA^n_0=(D_n,\emptyset)$ where \bae\label{d}
D_n=[-n,n]\cap \BZ\times\{0\}.
\eae 

Next, recall in \cite{procaccia2019stationary}, the stationary harmonic measure $\sharm$ on $\BH$ was defined as: for any $B\subseteq\BH$, any edge $\vec e=x\rightarrow y\in\partial^e B$, and any $N$,
$$
\sharm_{B,N}(\vec e)=\sum\limits_{z\in \ell_N\setminus B}\p_z(S_{\tau_{B\cup \ell_0}}=y,S_{\tau_{B\cup \ell_0}-1}=x).
$$
\begin{proposition}[Proposition 1, \cite{1}]
For any $B$ and $\vec e$ as above, there is a finite $\sharm_B(\vec e)$ such that
$$
\lim\limits_{N\to\infty}\sharm_{B,N}(\vec e)=\sharm_B(\vec e).
$$
\end{proposition}
 $\sharm_B(\vec e)$ is called the {\bf stationary harmonic measure} of $\vec e$ with respect to $B$ and the limit $\sharm_B(x)$ is called the {\bf stationary harmonic measure} of $x$ with respect to $B$. Then we give an informal description of the infinite SDLA model (see \cite{procaccia2019stationary} for details). Let $SV_0^\infty=\ell_0,SE_0^\infty=\emptyset$, and for any $t>0$, each edge $\vec e$ on the boundary of $SV_{t^-}^\infty$ is added to the edge set $SE_{t^-}^\infty$ and at the same time $\vec e(1)$ is added to the vertex set $SV_{t^-}^\infty$ at rate $\sharm_{SV_{t^-}^\infty}(\vec e)$. The process $SA_t^\infty=\left(SV_t^\infty,SE_t^\infty\right)$ starting from $\ell_0$ is called the infinite SDLA process. The following proposition says that $SA_t^\infty$ is well-defined. 
\begin{proposition}[Theorem 1, \cite{procaccia2019stationary}]
The infinite SDLA  $\{SA_t^\infty\}_{t\ge 0}$ is well defined.
\end{proposition}
Notice that there is a one-to-one correspondence between the elements in $\{G: G\subseteq\vec\BL^2\}$ and $\{\eta_G:\eta_G\in\{0,1\}^{\vec\BL^2}\}$ since for any directed subgraph $G=(V,\vec E)\subseteq \vec L^2$, we can define 
$$\eta_G(x)=
\begin{cases}
1& x\text{}\in G\\
0&\text{otherwise}
\end{cases},\
\eta_G(\vec e)=
\begin{cases}
1& \vec e\text{}\in G\\
0&\text{otherwise}
\end{cases}\ \forall (x,\vec e)\in G.
$$
So that both of the EDLA and SDLA process form Feller processes with sample paths in
$$D_E[0,\infty)=\left\{\text{right continuous functions $x: [0, \infty) \rightarrow E $ with left limits} \right\}$$ where $E=\{0,1\}^{\vec \BL^2}$.  
 The metric $\rho$ (defined in Section 4.1. of \cite{liggett2010continuous}) on $E$ induces a metric $d$ which gives rise to the Skorohod Topology on $D_E[0,\infty)$ (see Section 3.5 of \cite{ethier2009markov} for details). We say $\{EA^n_{nt}\cap \vec\BH\}_{t\ge 0}$ converges weakly to $\{SA^{\infty}_{ct}\}_{t\ge 0}$ iff their corresponding distributions converge.

With Remark \ref{remk1_1}, it is clear that the following theorem is an answer to Conjecture 1 of \cite{procaccia2019stationary}.
\begin{theorem}
\label{thm_scaling_limit}
	There exists $c\in (0,\infty)$ such that $EA^n_{nt}\cap \vec\BH$ converges weakly to $SA^{\infty}_{ct}$ on $(D_E[0,\infty),d)$ as $n\to\infty$, where $(D_E[0,\infty),d)$ is the metric space with the Skorohod topology. %\note{Don't we need to specify the topology of the weak convergence: weakly on compact time intervals and compact subsets of $\BH$}
\end{theorem}

\begin{notation}
	In this paper we will use $c,C$ etc. to denote constants. However, their values may vary according to contexts. 
\end{notation}

\begin{remark}
	The arguments in this paper also prove that the scaling limit of the regular DLA starting from a long line segment forms a variant of SDLA from $\ell_0$ where the growth rate is according to the stationary harmonic measure $\sharm$ on the outer boundary of the current aggregation.  
\end{remark}

\begin{remark}
	The SDLA or as shown in this paper the bulk of DLA stating from a long line, is expected to have a different fractal dimension from the standard DLA starting at a point. We conjecture that the dimension is $1.5$. This conjecture is based on connections to a stationary version of the Hastings Levitov process which is expected to have the same dimension. 
\end{remark}

It is easy to show the equivalence between the weak convergence and the finite dimensional distribution's convergence. So we put the proof of the following lemma in Appendix \ref{appendix}. 
\begin{lemma}\label{lem1_1}$EA^n_{nt}\cap \vec\BH$ converges weakly to $SA^{\infty}_{ct}$ if and only if the finite dimensional distribution of $EA^n_{nt}\cap \vec\BH$ converges to the corresponding finite dimensional distribution of $SA^{\infty}_{ct}$. Equivalently,
for any $\ep>0$, any finite subgraph $K\subseteq \vec\BH$ and $T<\infty$, there exists $N_0<\infty$ such that for any integer $n\ge 1$, $0<t_1,t_2,\cdots, t_n\le T$ and subgraph(s) $K_1,K_2,\cdots, K_n\subseteq K$, 
\bae\label{1_1}
\begin{aligned}
&\left|\prob\left(SA^\infty_{ct_1}\cap K=K_1, SA^\infty_{ct_2}\cap K=K_2,\cdots, SA^\infty_{ct_n}\cap K=K_n \right)\right.
\\& \hspace{0.2 in}\left.-\prob\left(EA^N_{Nt_1}\cap K=K_1, EA^N_{Nt_2}\cap K=K_2,\cdots, EA^N_{Nt_n}\cap K=K_n \right) \right|<\ep
\end{aligned}
\eae
for all $N\ge N_0$. 

\end{lemma}

  Let $SA^{m}_t$ be the SDLA process starting from $D_m$.
First by Theorem 1 of \cite{procaccia2019stationary}, $\{SA^{m}_t\}_{m\ge 1}$ and $SA^{\infty}_t$ can be coupled in the same probability space such that for any compact $K\subseteq \BH$ and any $T<\infty$, we have almost surely
\bae
\label{SDLA_converge}
SA^{m}_t\cap K\equiv SA^{\infty}_t\cap K, \ \forall t\in [0,T]
\eae
for all sufficiently large $m$. Thus in order to prove Theorem \ref{thm_scaling_limit}, by Lemma \ref{lem1_1}, it suffices to replace $SA^{\infty}_t$ with $SA^m_t$ and show the following proposition:

\begin{proposition}
\label{prop_limit}
For any $\ep>0$, any finite subgraph $K\subseteq \vec\BH$ and $T<\infty$, there exist $m_0,N_0<\infty$ such that for any integer $n\ge 1$, $0<t_1,t_2,\cdots, t_n\le T$ and subgraph(s) $K_1,K_2,\cdots, K_n\subseteq K$, 
\bae
\label{eq_close_1}
\begin{aligned}
&\left|\prob\left(SA^m_{ct_1}\cap K=K_1,  SA^m_{ct_2}\cap K=K_2,\cdots, SA^m_{ct_n}\cap K=K_n  \right)\right.\\
& \hspace{0.2 in}\left.-\prob\left(EA^N_{Nt_1}\cap K=K_1, EA^N_{Nt_2}\cap K=K_2,\cdots, EA^N_{Nt_n}\cap K=K_n  \right) \right|<\ep
\end{aligned}
\eae
for all $m\ge m_0$ and $N\ge N_0$. 
\end{proposition}

\subsection{The intermediate DLA process}
For the proof of Proposition \ref{prop_limit}, we introduce a family of intermediate DLA processes $IA^{m,N}_t$ defined as follows:
\begin{definition}
	\label{def_intermediate_process}
	For all positive integers $m\le N$, define the intermediate DLA process $IA^{m,N}_t=\left(IV^{m,N}_t,IE^{m,N}_t\right)$ to be a continuous time Markov process on the set of all subgraphs of $\vec\BL^2$ such that 
	\begin{itemize}
		\item $(IV^{m,N}_0,IE^{m,N}_0)=(D_m,\emptyset)$. 
		\item Assume there is a Poisson clock with intensity $N$. For any $s\ge 0$, if the clock rings at time $s$, we add $\vec e$ to $IE^{m,N}_{s-}$ and $\vec e(1)$ to $IV^{m,N}_{s-}$ such that 
		$$
		\left(IV^{m,N}_{s},IE^{m,N}_{s}\right)=\left(IV^{m,N}_{s-}\cup\{\vec e(1)\},IE^{m,N}_{s-}\cup\{\vec e\}\right)
		$$ with probability
				$$
		\eharm_{IV^{m,N}_{s-}\cup D_N}\left(\vec e\right) 		$$
	$\text {for all edges }\vec e\in \partial^{e} (IV^{m,N}_{s-})$.

	\end{itemize}
\end{definition}

It is clear that $IA^{m,N}_{t}$ forms a well defined (lazy) Markov process where a new particle is added at a rate uniformly bounded from above by $N$. 

First by a maximal coupling, we show that when $m,N$ is sufficiently large, $IA^{m,N}_t$ is the same as $SA^m_t$ with very high probability. That is,
\begin{proposition}
	\label{prop_limit_2}
	
 There exists $c>0$ such that for any $\ep>0,T<\infty,$ there is a constant $M_0<\infty$. And for all $m>M_0$ there exists $N(m)<\infty$ such that for all $N>N(m)$ we can couple $IA_{t}^{m,N}$ and $SA_t^m$ such that
 \bae\label{ep}
 P(IA_{t}^{m,N}\equiv SA_{ct}^m, \forall t\le T)\ge 1-\ep.
 \eae

\end{proposition}

Next, by coupling pairs of the intermediate DLA processes, we show that  for all $m\le N^{1/5}$, with high probability,  $IA_t^{m,N}$ and $IA_t^{m+1,N}$ have no discrepancy in $K$, when $m,N$ is sufficiently large. To be noted, $N^{1/5}$ is an adequate but not the only scale we can choose. 
\begin{proposition}
	\label{prop_coupling}
	For any finite subgraph $K\subseteq \vec\BH,T<\infty$, there exist $C<\infty$ and $\alpha>0$ such that for all sufficiently large $N,m$ satisfying $0<m\le N^{1/5}$, $IA^{m,N}_t$ and $IA^{m+1,N}_t$ can be coupled so that 
	\bae
	\label{eq_coupling_1}
		\prob\left(IA^{m,N}_t\cap K\equiv  IA^{m+1,N}_t\cap K, \forall t\le T\right)\ge 1-\frac{C}{m^{1+\alpha}}.
	\eae
\end{proposition}
 When $N$ is large enough, although $IA_t^{N^{1/5},N}$ and $IA_t^{N,N}$ seem to behave significantly differently near the end of the interval $D_N$, we can show that they are highly likely to be the same when restricted in a finite graph $K$. I.e.,
 \begin{proposition}
\label{prop_coupling2}
For any finite subgraph $K\subseteq \vec\BH$, any $\ep>0,T<\infty$, there exists $N_0>0$ such that for all $N\ge N_0$, $IA^{N^{1/5},N}_t$ and $IA^{N,N}_t$ can be coupled so that 
	\bae
	\label{eq_coupling_2}
		\prob\left(IA^{N^{1/5},N}_t\cap K\equiv  IA^{N,N}_t\cap K, \forall t\le T\right)\ge 1-\ep.
	\eae
\end{proposition}
\begin{notation}
Without loss of generality, we take $T=1$ in the rest of this paper.
\end{notation}

\subsection{Ideas and structure of the proof }
At first, we explain how to establish Proposition \ref{prop_limit} from Proposition \ref{prop_limit_2}-\ref{prop_coupling2}. Fix a sufficiently large $m$, a finite graph $K$, it is sufficient for us to find $N_0$ such that \eqref{eq_close_1} holds for all $N\ge N_0$. Proposition \ref{prop_limit_2} tells us that  there exists $N_m$ such that for all $N\ge N_m$, $IA_t^{m,N}\equiv SA_{ct}^m$ on $[0,1]$  with high probability. Then Proposition \ref{prop_coupling} tells us that we can find $\tilde N_m\ge \max\{m^5,N_m\}$ such that for all $N\ge \tilde N_m$,  with small probability there exists $m\le \tilde m\le N^{1/5}$ such that $IA_t^{\tilde m,N}\cap K\not\equiv IA_t^{\tilde m+1,N}\cap K$ on $[0,1]$. At last, Proposition \ref{prop_coupling2} tells us that we can find $\hat N_m>\tilde N_m$ such that for all $N\ge \hat N_m$, $IA_t^{N^{1/5},N}\cap K\equiv IA_t^{N,N}\cap K$ on $[0,1]$ with high probability. Then we can choose $N_0=\hat N_m$. 

  To couple all the finite discrete intermediate DLA processes $\{IA^{m,N}_k\}_{k\le 2N},m\le N$ together, we sample $2N$ i.i.d. copies of SRW's starting from the outer boundary of the ball $B(0,4N)$ according to the regular harmonic measure $\CH$ and accomplish the task in Section \ref{sec3}.

In Section \ref{sec_upper bound}, we obtain upper bounds on the growth of the intermediate DLA processes. As a result, we only need to consider the truncated processes without growing outside a finite region in the following sections.

We begin to prove our result in Section \ref{sec_prop_limit_2}. First we show Proposition \ref{prop_limit_2}. There we consider the truncated continuous time coupled process $(IA_t^{m,N},SA_t^m)$ constructed by a maximal coupling. By Lemma \ref{lem1}, when $IA_{(t\land\Gamma_m)-}^{m,N}=SA_{(t\land\Gamma_m)-}^m$, the total transition rate of $(IA_{t\land\Gamma_m}^{m,N},SA_{t\land\Gamma_m}^m)$ converges to 0 uniformly in the unit time interval. Since $IA_{0}^{m,N}=SA_{0}^m$, we obtain that the probability $IA_{t\land\Gamma_m}^{m,N}\equiv SA_{t\land\Gamma_m}^m$ on $[0,1]$ converges to 0 when $m,N$ converges to infinity.  

In the last two sections, Section \ref{sec_prop_coupling} and \ref{Sec_prop_coupling2}, we consider the discrete time truncated coupled process $(IA_{k\land\Gamma_m}^{m,N},IA_{k\land\Gamma_m}^{m+1,N})$ and prove Proposition  \ref{prop_coupling} and \ref{prop_coupling2}. The idea of those two sections borrows techniques from \cite{procaccia2019stationary}, which concentrated on the continuous time process. We trace the positions of the two edge discrepancies $\vec e_{\Delta_i,1}, \vec e_{\Delta_i,2}$ created at time $\Delta_i$, and show that in the $2N$ steps, the discrepancies do not reach any finite graph $K$ with high probability.

% !TEX root = Convergences_to_SDLA_3.12.19.tex

\section{Coupling construction}
\label{sec3}
Given $N$, let $IA^{m;N}_0=(D_m,\emptyset)$ for all $m\le N$. Let $\left\{S^{(k)}_n\right\}_{n=0}^\infty$,$1\le k\le 2N$ be $2N$ i.i.d. copies of SRW's starting at radius $4N$ according to the regular harmonic measure $\CH$. Then for any $1\le k\le 2N$, let $\tau^{(k)}$ be the stopping time with respect to $S^{(k)}$. 
\begin{itemize}
\item If 
$$
\tau^{(k)}_{IA^{m;N}_{k-1}}<\tau^{(k)}_{D_N\setminus D_m},
$$
we add the directed edge $S^{(k)}_{\tau^{(k)}_{IA^{m;N}_{k-1}}-1}\to S^{(k)}_{\tau^{(k)}_{IA^{m;N}_{k-1}}}$ to the edge set $IE_k^{m;N}$ and vertex $S^{(k)}_{\tau^{(k)}_{IA^{m;N}_{k-1}}-1}$ to the vertex set $IV_k^{m;N}$.
\item Otherwise, we keep $IA_k^{m;N}$ the same. 
\end{itemize}

So now we have coupled all $\{IA_k^{m;N}\}_{0\le k\le 2N}$, $m\le N$ together. By definition, for each $m\le N$, the marginal distribution of $IA_k^{m;N}$ is the embedded chain of the intermediate DLA process.

\begin{remark}\label{rem4_2}
 By large deviation principle, with high probability the transitions for $IA_{t}^{m,N}$ in the unit time is no more than $2N$ since the waiting time of each transition has the exponential distribution $\exp(N)$. That's why we consider the finite embedded chain $IA_{k}^{m,N},k\le 2N$. \end{remark}

Now we concentrate on the distribution of the pair $(IA^{m_1;N},IA^{m_2;N}),m_1<m_2$ which plays an important role in the proofs of Proposition \ref{prop_coupling} and \ref{prop_coupling2}. 
Define
 $$
 \eharm_A\left(x,\vec e\right)=\p_x\left(\bar\tau_A=\bar\tau_{\vec e\left(2\right)},S_{\bar\tau_{\vec e\left(2\right)}-1}=\vec e\left(1\right)\right)
 $$
 and for any subgraph $G=(V,E)\subseteq\vec \BL^2$, and any directed edge $\vec e\in \vec\BL^2$, denote
 $$
 G\cup\{\vec e\}=\left(V\cup\{\vec e(1),\vec e(2)\},E\cup\{\vec e\}\right).
 $$
Formally, the construction of the coupled Markov chain $(IA^{m_1;N},IA^{m_2;N}),k\le 2N$ is described as follows:

 \begin{itemize}
\item $\left(IA_0^{m_1,N},IA_0^{m_2,N}\right)=\left(\left(D_{m_1},\emptyset\right),\left(D_{m_2},\emptyset\right)\right)$.
\item For any $1\le k\le 2N$, denote the joint transition probability that from $\left(IA_k^{m_1,N},IA_k^{m_2,N}\right)$ to $\left(IA_{k+1}^{m_1,N},IA_{k+1}^{m_2,N}\right)$ as 
$$\p\left(\left(IA_k^{m_1,N},IA_k^{m_2,N}\right),\left(IA_{k+1}^{m_1,N},IA_{k+1}^{m_2,N}\right)\right).$$ 
\end{itemize}
Then if they exist, we define the first added edge at time $k$ as $\vec e_{k,1}$ and the second added edge as $\vec e_{k,2}$, so that
$$
\vec e_{k,i}=S^{(k)}_{\tau^{(k)}_{IA^{m_i;N}_{k-1}}-1}\to S^{(k)}_{\tau^{(k)}_{IA^{m_i;N}_{k-1}}},
\ i=1,2.$$
  Then there are eight cases that may happen. In the first three cases, there are two added edges added at time $k$, while in the rest five cases, $S_n^{(k)}$ hits $D_N$ before the second edge is added so that there is at most one edge added. Especially, in the last case, $S_n^{(k)}$ hits $D_N$ before $\min\left\{\tau^{(k)}_{IA^{m_1;N}_{k-1}},\tau^{(k)}_{IA^{m_2;N}_{k-1}}\right\}$, so that no edge is added.
   \begin{enumerate}[I.]
 
\item  If $\vec e_{k,1}\left(2\right)\in IA_k^{m_1,N}\cap IA_k^{m_2,N} $, 
 we have$$ 
\left(IA_{k+1}^{m_1,N},IA_{k+1}^{m_2,N}\right)=\left(IA_k^{m_1,N}\cup\left\{\vec e_{k,1}\right\},IA_k^{m_2,N}\cup\left\{\vec e_{k,1}\right\}\right)
$$
and
$$
\p\left(\left(IA_k^{m_1,N},IA_k^{m_2,N}\right),\left(IA_{k+1}^{m_1,N},IA_{k+1}^{m_2,N}\right)\right)=\eharm_{IV_k^{m_1,N}\cup IV_k^{m_2,N}\cup D_N}\left(\vec e_{k,1}\right).
$$
\item  If $\vec e_{k,1}\left(2\right)\in IA_k^{m_1,N}\cap\left( IA_k^{m_2,N} \cup D_N\right)^c, \vec e_{k,2}(2)\in IV_k^{m_2,N}$, 
we have $$
\left(IA_{k+1}^{m_1,N},IA_{k+1}^{m_2,N}\right)=\left(IA_k^{m_1,N}\cup\left\{\vec e_{k,1}\right\},IA_k^{m_2,N}\cup\left\{\vec e_{k,2}\right\}\right)
$$and
\bae\label{4_1}
\p\left(\left(IA_k^{m_1,N},IA_k^{m_2,N}\right),\left(IA_{k+1}^{m_1,N},IA_{k+1}^{m_2,N}\right)\right)=\eharm_{IV_k^{m_1,N}\cup IV_k^{m_2,N}\cup D_N}\left(\vec e_{k,1}\right)\eharm_{IV_k^{m_2,N}\cup  D_N}\left(\vec e_{k,1}\left(2\right),\vec e_{k,2}\right).
\eae
\item  If $\vec e_{k,1}\left(2\right)\in IA_k^{m_2,N}\cap\left( IA_k^{m_1,N}\cup D_N \right)^c, \vec e_{k,2}(2)\in IA_k^{m_1,N}$, 
we have $$
\left(IA_{k+1}^{m_1,N},IA_{k+1}^{m_2,N}\right)=\left(IA_k^{m_1,N}\cup\left\{\vec e_{k,2}\right\},IA_k^{m_2,N}\cup\left\{\vec e_{k,1}\right\}\right)
$$and
\bae\label{4_2}
\p\left(\left(IA_k^{m_1,N},IA_k^{m_2,N}\right),\left(IA_{k+1}^{m_1,N},IA_{k+1}^{m_2,N}\right)\right)=\eharm_{IV_k^{m_1,N}\cup IV_k^{m_2,N}\cup D_N}\left(\vec e_{k,1}\right)\eharm_{IV_k^{m_1,N}\cup  D_N}\left(\vec e_{k,1}\left(2\right),\vec e_{k,2}\right).
\eae

\item If $\vec e_{k,1}\left(2\right)\in IA_k^{m_1,N}\cap\left( IA_k^{m_2,N} \cup D_N\right)^c, \vec e_{k,2}\left(2\right)\in D_N\backslash  D_{m_2}$, 
we have $$
\left(IA_{k+1}^{m_1,N},IA_{k+1}^{m_2,N}\right)=\left(IA_k^{m_1,N}\cup\left\{\vec e_{k,1}\right\},IA_k^{m_2,N}\right)
$$and
\bae\label{4_9}
\p\left(\left(IA_k^{m_1,N},IA_k^{m_2,N}\right),\left(IA_{k+1}^{m_1,N},IA_{k+1}^{m_2,N}\right)\right)=\eharm_{IV_k^{m_1,N}\cup IV_k^{m_2,N}\cup D_N}\left(\vec e_{k,1}\right)\eharm_{IV_k^{m_2,N}\cup  D_N}\left(\vec e_{k,1}\left(2\right),\vec e_{k,2}\right).
\eae
\item  If $\vec e_{k,1}\left(2\right)\in IA_k^{m_2,N}\cap\left( IA_k^{m_1,N}\cup D_N \right)^c, \vec e_{k,2}\left(2\right)\in D_N\backslash  D_{m_1}$, 
we have $$
\left(IA_{k+1}^{m_1,N},IA_{k+1}^{m_2,N}\right)=\left(IA_k^{m_1,N},IA_k^{m_2,N}\cup\left\{\vec e_{k,1}\right\}\right)
$$and
\bae\label{4_10}
\p\left(\left(IA_k^{m_1,N},IA_k^{m_2,N}\right),\left(IA_{k+1}^{m_1,N},IA_{k+1}^{m_2,N}\right)\right)=\eharm_{IV_k^{m_1,N}\cup IV_k^{m_2,N}\cup D_N}\left(\vec e_{k,1}\right)\eharm_{IV_k^{m_1,N}\cup  D_N}\left(\vec e_{k,1}\left(2\right),\vec e_{k,2}\right).
\eae
\item  If $\vec e_{k,1}\left(2\right)\in IA_k^{m_1,N}\cap\left( IA_k^{m_2,N} \right)^c\cap D_N$, 
we have $$
\left(IA_{k+1}^{m_1,N},IA_{k+1}^{m_2,N}\right)=\left(IA_k^{m_1,N}\cup\left\{\vec e_{k,1}\right\},IA_k^{m_2,N}\right)
$$and
\bae\label{4_11}
\p\left(\left(IA_k^{m_1,N},IA_k^{m_2,N}\right),\left(IA_{k+1}^{m_1,N},IA_{k+1}^{m_2,N}\right)\right)=\eharm_{IV_k^{m_1,N}\cup IV_k^{m_2,N}\cup D_N}\left(\vec e_{k,1}\right).
\eae
\item  If $\vec e_{k,1}\left(2\right)\in IA_k^{m_2,N}\cap\left( IA_k^{m_1,N} \right)^c\cap D_N$, 
we have $$
\left(IA_{k+1}^{m_1,N},IA_{k+1}^{m_2,N}\right)=\left(IA_k^{m_1,N},IA_k^{m_2,N}\cup\left\{\vec e_{k,1}\right\}\right)
$$and
\bae\label{4_12}
\p\left(\left(IA_k^{m_1,N},IA_k^{m_2,N}\right),\left(IA_{k+1}^{m_1,N},IA_{k+1}^{m_2,N}\right)\right)=\eharm_{IV_k^{m_1,N}\cup IV_k^{m_2,N}\cup D_N}\left(\vec e_{k,1}\right).
\eae
\item  Otherwise, we have
$$\left(IA_{k+1}^{m_1,N},IA_{k+1}^{m_2,N}\right)=\left(IA_k^{m_1,N},IA_k^{m_2,N}\right).
$$
 
\end{enumerate}

 Now we use the definition of the vertex discrepancies and edge discrepancies in \cite{procaccia2019stationary} such that
$$
V_n^{D,m_1,m_2}=\left\{x\in\BZ^2: s.t. \ \exists k\le n,x\in  IV_k^{m_1,N}\triangle IV_k^{m_2,N} \right\}
$$
denotes the the set of vertex discrepancies and 
$$
E_n^{D,m_1,m_2}=\left\{\vec{e}\in\BZ^2: s.t. \ \exists k\le n,\vec e\in  IE_k^{m_1,N}\triangle IE_k^{m_2,N} \right\}
$$
denotes the set of edge discrepancies before time $n$ where $\triangle$ stands for the symmetric difference between sets.
  From the definition above, we give the following statement to deepen our understanding on their relations.
\begin{itemize}
\item For any vertex $x\in V_n^{D}\cap A$, there must be an edge $\vec e$ in $  E_n^{D}\cap \left(\vec A\cup\widetilde{\partial^e A}\right)$ such that $x=\vec e(1)$.
\item For any $\vec e$ in $E_n^{D}\cap \left(\vec A\cup\widetilde{\partial^e A}\right)$, $\vec e(1)\in V_n^{D}\cap A$.
\end{itemize}

Denote the stopping times enumerating discrepancies as 
  \bae
 & \Delta^{m_1,m_2}_1=\inf\left\{1\le k\le 2N:|E_k^{D,m1,m_2}\backslash E_{k-1}^{D,m_1,m_2}|\ge 1\right\},
 \\& \Delta^{m_1,m_2}_i=\inf\left\{ 2N\ge k> \Delta^{m_1,m_2}_{i-1}:|E_k^{D,m_1,m_2}\backslash E_{k-1}^{D,m_1,m_2}|\ge 1\right\}
 \eae
and with convention that $\inf\emptyset=\infty$.
 Denote the set of all the stopping times as $T_\Delta^{m_1,m_2}$. 
 
\begin{remark}\label{rem4_1}
 Note that the event $\left\{n\in T_\Delta^{m_1,m_2}\right\}$ is equivalent to the event $$\left\{\vec e_{n,1}(2)\in IV^{m_1,N}_{n-1}\triangle IV^{m_2,N}_{n-1}\subseteq V_{n-1}^{D,m_1,m_2}\right\},$$ whose probability is the summation over probabilities represented in $\eqref{4_1}$-$ \eqref{4_10}$.
 \end{remark}

% !TEX root = Convergences_to_SDLA_3.12.19.tex
\section{Upper bounds on the growth of the intermediate processes}
\label{sec_upper bound}
Before proving our results, we first give some useful lemmas, mainly the upper bounds on the edge harmonic measure Lemma \ref{lem2} and the growth rates of the intermediate DLA processes, Lemmas \ref{lem2_1} and \ref{lem2_2}. Given these estimates, we will only need to consider a truncated processes in a finite region. 

The first lemma is about the stochastic domination of independent Bernoulli random variables. It is very simple to prove by induction, whence one who has interests can refer to Appendix \ref{appendix}.
\begin{lemma}\label{lem2-1}
If $X_1,\cdots,X_n$ are $n$ random variables satisfying that 
$$
\p\left(X_1=1\right)\le p,\text{ }\p\left(X_k=1|X_1=a_1,\cdots,X_{k-1}=a_{k-1}\right)\le p
$$
for any $\left(a_1,\cdots,a_{k-1}\right)\in\{0,1\}^{k-1},2\le k\le n$, then $X_1,\cdots,X_n$ can be stochastically dominated by independent Bernoulli random variables $Y_1,\cdots,Y_n$ with parameter $p$.
\end{lemma}

Denote
$$
F_m=[-m-\log m,m+\log m]\times [-\log m,\log m]\cap\BZ^2.
$$

Next we give an upper bound on the rescaled edge harmonic measure $N\eharm_{A\cup D_N}\left(y\right)$ for all $y$ in a thin subset $ F_m$. Since the proof of Lemma \ref{lem2} is very similar to existing results from the literature we also push it to Appendix \ref{appendix}. 
\begin{lemma}\label{lem2}
For any $\delta>0$, $m\le\left(1-\delta\right)N$, and $x\in F_m$, there exists $C\in\left(0,\infty\right)$ which is independent of $A$ such that for any connected $A\subseteq\BZ^2$ with $D_N\subseteq A$,
$$
N \eharm_{A\cup D_N}\left(x\right)\le C\sqrt{|x\left(2\right)|}
$$ when $N$ is sufficiently large.
\end{lemma}
We will make use of a uniform upper bound on the regular harmonic measure proved by Kesten in 1987. 
\begin{lemma}[Theorem of \cite{kesten1987hitting}]
\label{citelem}
Let $A$ be a connected subset in $\BZ^d$ which contains the origin. Then there exists a constant $C_0\in \left(0,\infty\right)$, independent of $A$, such that for all $x\in A$,
$$
\CH_A\left(x\right)\le C_0||A||^{-1/2}
$$where $\|A\|$ is the radius of $A$.
\end{lemma}
Define two boxes 
\bae
&B_1=\left([-N-4C_0N^{1/2},N/2]\cup [N/2,N+4C_0N^{1/2}]\right)\times[-4C_0N^{1/2}, 4C_0N^{1/2}]\cap \BZ^2,
\\&B_2=[-N/2,N/2]\times[-\log N, \log N]\cap\BZ^2.
\eae
Next we will explain how the upper bound on the growth rate fit in proving the logarithm growth upper bound for the intermediate process with a long boundary.
\begin{lemma}\label{lem2_1}For any $C_1<\infty,\delta>0$, $m\le \left(1-\delta\right)N$,

$$\p\left(IA_{2N}^{m,N}\subseteq \vec F_m\right)>1-\frac{1}{m^{C_1}}
$$for all sufficiently large $N$.\end{lemma}
\begin{proof}
Denote $IA^{m,N}_k\left(x\right)$ as the connected component of $x$ in $IA^{m,N}_k$ such that its vertex set
$$
IV^{m,N}_k\left(x\right)=\left\{y\in\BH: x\text{ is connected to }y \text{ by a directed path in }IA^{m,N}_k\right\}.
$$ Then it is easy to see that
\bae\label{2_00}IA_{2N}^{m,N}=\cup_{x\in D_m}IA_{2N}^{m,N}\left(x\right).\eae
For any $x\in D_m$, if $IV_{2N}^{m,N}\left(x\right)\cap F_m^c\ne\emptyset$, there must be a nearest neighbor directed path in $IA_{2N}^{m,N}\left(x\right)$ such that 
$$
\CP_x=\{P_{\log m}\to P_{\log m-1}\to\cdots\to P_0=x\},||P_i-P_{i-1}||=1,0<i\le \log m
$$ from some point $P_{\log m}$ in $F_m$ to $x$. 
Define the random variable
$$
X_n=
\begin{cases}
1&\text{ if $P_i\in IV_n^{m,N}\left(x\right)$ for some $1\le i\le \log m$ and $\CP_x\cap IV_{n-1}^{m,N}\left(x\right)=\{P_0,\cdots,P_{i-1}\}$}\\
0&\text{ otherwise}
\end{cases}
$$forall $1\le n\le 2N$.

By Lemma \ref{lem2}, 
$$
\p[X_n=1|\CF_{n-1}]\le\frac{C\sqrt{\log m}}{N}
$$
where $\CF_n$ is the $\sigma$-field generated by $IA_k^{m,N},k\le n$.
So that by Lemma \ref{lem2-1}, $\{X_n,1\le n\le 2N\}$ can be stochastically dominated by the independent  random variables $\{Y_n,1\le n\le 2N\}$ which satisfies
$$
\p\left(Y_n=1\right)=1-\p\left(Y_n=0\right)=\frac{C\sqrt{\log m}}{N}.
$$ 
It follows that for any $\theta>0$
\bae\label{2_10}
\p\left(\sum\limits_{n=1}^{2N}X_n\ge \log m\right)
&\le\p\left(\sum\limits_{n=1}^{2N}Y_n\ge \log m\right)
\\&\le \frac{\e \exp\left(\theta\sum\limits_{n=1}^{2N}Y_n\right)}{\exp\left(\theta\log m\right)}
\\&=\frac{\left(1+ C[\exp\left(\theta\right)-1]\sqrt{\log m}/N\right)^{2N}}{\exp\left(\theta\log m\right)}
\\&\sim \exp\left(C\left(\theta\right)\sqrt{\log m}-\theta\log m\right)
\eae
 when $N$ is large enough where $C\left(\theta\right)$ is a constant associated with $\theta$.
By \eqref{2_00} and \eqref{2_10}, for any $C_1<\infty$, 
\bae
\p\left(IV_{2N}^{m,N}\not\subseteq F_m\right)&=\p\left(\cup_{x\in D_m}IV_{2N}^{m,N}\left(x\right)\not\subseteq F_m\right) \\&\le2m\p\left(\CP_0\text{ exists }\right)
\\&\le 2m4^{\log m}\exp\left(C\left(\theta\right)\sqrt{\log m}-\theta\log m\right)
\\&\le\exp\left(-C_1\log m\right)
\eae
  when $m$ is large enough, where the last inequality holds by choosing an adequate $\theta$.
\end{proof}
The next lemma gives an upper bound on the probability that the sum of uniformly bounded independent random variables deviates from its conditional expectations given the past. It will be used plenty of times in the following proofs.
\begin{lemma}[Theorem of \cite{freedman1973some}]\label{lem2_1,1}
Suppose $0\le X_i\le 1$ and $X_i $ is $\CF_i$ measurable. Let $M_i=\e\left(X_i|\CF_{i-1}\right)$, for any $0\le b\le a$
$$\p\left(\sum\limits_{i=1}^{n}X_i\ge a,\sum\limits_{i=1}^{n}M_i\le b\right)\le\exp\left(-\frac{\left(a-b\right)^2}{2a}\right).$$
\end{lemma}
Note that the logarithm growth does not hold when $m=N$, i.e. $IA_{t}^{N,N}=EA_{Nt}^{N}$. But we can still give a rough upper bound on the growth of $IA_{2N}^{N,N}$ which is good enough for our proof.
\begin{lemma}\label{lem2_2}
For any $C<\infty$,

$$\p\left( IA_{2N}^{N,N}\subseteq \vec B_1\cup \vec B_2\right)>1-\frac{1}{N^C}
$$for all sufficiently large $N$.\end{lemma}

\begin{proof}
Similar to Lemma \ref{lem2_1}, we can prove that for any $C_1\in \left(0,\infty\right)$, 
 \bae\label{2_9}
 \p\left(\cup_{x\in D_{2N/3}}IV_{2N}^{N,N}\left(x\right)\subseteq B_1\cup B_2\right)\ge 1-\frac{1}{N^{C_1}}.
 \eae
 Thus conditional on the event $\left\{\cup_{x\in D_{2N/3}}IV_{2N}^{N,N}\left(x\right)\subseteq B_1\cup B_2\right\}$, if $IV_{2N}^{N,N}\cap\left( B_1\cup B_2\right)^c\ne\emptyset$, we must have 
 $$
\left(\cup_{x\in D_N\backslash D_{2/3N}} IV_{2N}^{N,N}\left(x\right)\right)\cap\left( B_1\cup B_2\right)^c\ne\emptyset.
$$
So that there must be a nearest neighbor directed path in $IA_{2N}^{N,N}\left(x\right)$ with $x\in D_N\backslash D_{2N/3}$ such that 
$$
\CP_x=\{P_{4C_0\sqrt{N}}\to P_{4C_0\sqrt{N}-1}\to\cdots\to P_0=x\},||P_i-P_{i+1}||=1,0\le i\le 4C_0\sqrt{N}
$$ from some point $P_{4C_0\sqrt{N}}$ in $B_1\cup B_2$ to $x$.

Define random variable
$$
X_n=
\begin{cases}
1&\text{ if $P_i\in IV_n^{N,N}\left(x\right)$ for some $1\le i\le 4C_0\sqrt{N}$ and $\CP_x\cap IV_{n-1}^{N,N}\left(x\right)=\{P_0,\cdots,P_{i-1}\}$}\\
0&\text{ otherwise}
\end{cases}
$$for all $1\le n\le 2N$.
 
By Lemma \ref{citelem}, $\forall 1\le n\le 2N$,
 $$
\p[X_n=1|\CF_{n-1}]\le \frac{C_0}{\sqrt{N}}.
$$
And by Lemma \ref{lem2_1,1},
\bae\label{2_1}
&\p\left(\#\{1\le n\le 2N:X_n=1\}\ge 4C_0\sqrt{N}\right)
\\&=\p\left(\#\{1\le n\le 2N:X_n=1\}\ge 4C_0\sqrt{N},\sum\limits_{n=1}^{2N}\p[X_n=1|\CF_{n-1}]\le C_0\sqrt{N}\right)
\\&\le \exp{\left(-C_0\sqrt{N}\right)}.
\eae
We deduce from \eqref{2_9} and \eqref{2_1} that for any $C<\infty$,
\bae
&\p\left(IA_{2N}^{N,N}\not\subseteq B_1\cup B_2 \right)\\&\le \p\left(\cup_{x\in  D_{2/3N}}IA_{2N}^{N,N}\left(x\right)\not\subseteq B_1\cup B_2\right) +\p\left(\cup_{x\in  D_N\backslash D_{2/3N}}IA_{2N}^{N,N}\left(x\right)\not\subseteq B_1\cup B_2\right)
\\&\le \frac{1}{N^{C_1}}+\frac{2N}{3}4^{C_0\sqrt{N}} \exp{\left(-C_0\sqrt{N}\right)}
\\&\le \frac{1}{N^{C}}
\eae
when $N$ is large enough.
\end{proof}

% !TEX root =Convergences_to_SDLA_3.12.19.tex

\section{Proof of Proposition \ref{prop_limit_2}}\label{sec_prop_limit_2}
In this section, we consider the continuous time process. First for completeness we state the following lemma, an adaption of Theorem 1.3 of \cite{procaccia2018stationary}.
\begin{lemma}[Adaption of Theorem 1.3 of \cite{procaccia2018stationary}]
\label{lem1}
For any finite connected subset $A\subseteq\BH$, there is a constant $C\in\left(0,\infty\right)$, independent of the set A, such that for any point $x\in A\backslash l_0$,
\bae
C\lim\limits_{n\to\infty}N\eharm_{A\cup D_N}\left(x\right)= \harm^s_{A\cup l_0}\left(x\right).
\eae 
Moreover, $C=2/\lim\limits_{n\to\infty}n\eharm_{D_n}\left(0\right)$.
\end{lemma}

Now we come to the main proof of this section.\\
 %{\bf Proof of Proposition \ref{prop_limit_2}:}\\
 \begin{proof}[Proof of Proposition \ref{prop_limit_2}]
Here we use the maximal coupling constructed in Section 1 of Chapter III of \cite {liggett1985}. Let $c=1/C$, where $C$ is the positive constant in Lemma \ref{lem1}. 
Define
$$
\Gamma_m=\inf\{t:IA_{t}^{m,N}\cup SA_{ct}^m\not\subseteq \vec F_m\}.
$$
For any $C\in\left(0,\infty\right)$, when $m$ is large enough, by Theorem 5 of \cite{procaccia2019stationary}, 
\bae\label{3_1}
\p\left(\exists t\le 1, SA_{ct}^m\not\subseteq  \vec F_m\right)\le \frac{1}{m^C},
\eae
while by Lemma \ref{lem2_1}, 
\bae\label{3_3}
\p\left( IA_{2N}^{m,N}\not\subseteq \vec F_m\right)\le \frac{1}{m^C}.
\eae

However, by the characteristic function of the Poisson distribution,
\bae\label{3_4}
&\p\left( \text{ there are more than $2N$ transitions up to time 1}\right)
\\&=\p\left(X\ge 2N\right)
\\&\le \frac{\e\exp(X)}{\exp\left(2N\right)}
\\&=\exp\left(-\left(3-e\right)N\right)
\eae
where $X$ is distributed Poisson(N).
We deduce from \eqref{3_1},\eqref{3_3} and \eqref{3_4} that for any $\ep>0$,
\bae\label{3_2}
\p\left(\Gamma_m\le 1\right)\le  \ep/2
\eae
when $m$ is large enough.

The truncated processes $IA_{t\land \Gamma_m}^{m,N}$ and $SA^m_{ct\land \Gamma_m}$ are two finite Markov processes  on $\{0,1\}^{ \vec F_m}$. We denote them as $\hat A_t^{m,N}$ and $\hat B^m_{ct}$ respectively. Considering the coupled process $Z_t=\left(\hat A_t^{m,N},\hat B^m_{ct}\right)$ on $\{0,1\}^{ \vec F_m}\times \{0,1\}^{ \vec F_m}$, by Lemma \ref{lem1} we have
\bae
&\lim\limits_{\Delta t\to0}\frac{\p\left(\exists s\le t+\Delta t, \hat A_s^{m,N}\ne \hat B^m_{cs}\right)-\p\left(\exists s\le t, \hat A_s^{m,N}\ne \hat B^m_{cs}\right)}{\Delta t}
\\&=\lim\limits_{\Delta t\to0}\frac{\p\left(\exists t< s\le t+\Delta t,\hat A_s^{m,N}\ne \hat B^m_{cs},\forall s\le t, \hat A_s^{m,N}\equiv\hat B^m_{cs}\right)}{\Delta t}
\\&\le\sup_{A\subseteq F_m}\sum\limits_{\vec{e}}|c\harm^s_{A\cup l_0}\left(\vec{e}\right)-N\eharm_{A\cup D_N}\left(\vec{e}\right)|
\\&\to 0
\eae
uniformly in $t\le 1$ when $N\to\infty$.

It follows that for any $\ep>0,m<\infty$, there exists $N_m$ such that for all $N\ge N_m$ we have 
\bae\label{20}
P\left(\hat A_t^{m,N}\not\equiv \hat B^m_{ct} \text{ on }[0,1]\right)\le \int_0^1\ep/2ds\le\ep/2.
\eae

Thus it follows from \eqref{3_2} and \eqref{20} that \eqref{ep} is true when $T=1$.
%$\hfill\square$ 
\end{proof}

% !TEX root = Convergences_to_SDLA_3.12.19.tex
\section{Proof of Proposition \ref{prop_coupling}}\label{sec_prop_coupling}
Recall the coupled process $$\left(IA_{k}^{m,N}, IA_{k}^{m+1,N}\right),k\le 2N$$ constructed in Section \ref{sec3}.
Define the stopping time
$$
\Gamma_m=\inf\left\{n\le 2N:IA_{n}^{m,N}\cup IA_{n}^{m+1,N}\not\subseteq  \vec F_{m+1}\right\},
$$
and the truncated process
$$
\left(\hat A^{m,N}_k, \hat A^{m+1,N}_k\right)=\left(IA_{k\land \Gamma_m}^{m,N}, IA_{k\land \Gamma_m}^{m+1,N}\right).
$$
By Lemma \ref{lem2_1}, for any $C\in \left(0,\infty\right)$ and sufficiently large $m$,
$$
\p\left(\Gamma_m<2N\right)<\frac{2}{m^C}.
$$
Then it suffices to show that for all sufficiently large $m$ satisfying $m\le N^{1/5}$, there exist $\alpha>0$ and $ C<\infty$ such that for any finite subgraph $K\subseteq \vec\BH$,
$$
\p\left(\exists k\le 2N,\hat A_k^{m,N}\cap K\ne\hat A_k^{m+1,N}\cap K\right)\le\frac{C}{m^{1+\alpha}}.
$$
Recall the definition of the stopping time $\Delta^{m_1,m_2}$ when a discrepancy occurs in Section \ref{sec3}. Let $T_\Delta^m$ be the set of the stopping times before $2N\land \Gamma_m$ and we abbreviate $\Delta^{m,m+1}_i$ to $\Delta^{}_i$ here, so that
$$
T_\Delta^m=\left\{\Delta^{}_i:\Delta^{}_i\le 2N\land \Gamma_m\right\}.
$$Then we want to get an upper bound on the number of the stopping times in $T_\Delta^m$.
\begin{lemma}\label{lem5_1}
For any $\alpha>0$, there exists $ c>0$ such that
\bae\label{eq_lem5_1}
\p\left(|T^m_\Delta|\ge m^\alpha\right)\le \exp\left(-m^c\right)
\eae for all sufficiently large $m,N$ with $m\le N^{1/5}$.
\end{lemma}
\begin{proof}

By Lemma \ref{lem2} and Remark \ref{rem4_1},
\bae\label{5_1}
\p\left(n\in T_\Delta^m|\CF_{n-1}\right)&\le\frac{C\sqrt{\log m}}{N}|V^{D,m}_{n-1}|
\eae
when $m,N$ are large enough and $m\le N^{1/5}$.
For any $\delta<1$, let $\Delta_0=0$ and
$$
\forall 1\le i\le m^\alpha,X_i=
\begin{cases}
1&\text{ if $\Delta_i-\Delta_{i-1}\le \frac{\delta N}{2i\sqrt{\log m}}$ or $\Delta_i=\infty$}\\
0&\text{ otherwise}
\end{cases}.
$$
Define
\begin{align*}
&I_k=\left\{\left(k-1\right)m^{\alpha/2}+1,\cdots,km^{\alpha/2}\right\},
\\&A_k=\left\{\sum\limits_{i\in I_k}X_i<c_0m^{\alpha/2}\right\},\forall 1\le k\le m^{\alpha/2}
\end{align*}for some $c_0>0$.
On $\cap_{1\le k\le m^{\alpha/2}}A_k$,
$$
\sum\limits_{i=1}^{m^\alpha}\Delta_i-\Delta_{i-1}\ge\sum\limits_{k=1}^{m^{\alpha/2}}c_0m^{\alpha/2} \times\frac{\delta N}{2km^{\alpha/2}\sqrt{\log m}}\ge\frac{c_0\delta\alpha N\sqrt{\log m}}{4}> 2N
$$for any $c_0,\delta>0$ when $m,N$ is sufficiently large enough.
It implies that
\bae\label{5_01}
\p\left(|T^m_\Delta|\ge m^\alpha\right)\le\p\left(\sum\limits_{i=1}^{m^\alpha}\Delta_i-\Delta_{i-1}\le 2N\right)\le\p\left(\cup_{1\le k\le m^{\alpha/2}}A_k^c\right).
\eae

Then it suffices to prove that for any $\alpha>0$, there exists $c>0$ such that
\bae\label{suffice}
\p\left(A_k^c\right)\le \exp\left(-m^c\right).
\eae
Notice that by strong Markov property,
\bae\label{5_1,1}
\p\left(X_i=1|\CF_{\Delta_{i-1}}\right)&=\p\left(\Delta_i-\Delta_{i-1}\le \frac{\delta N}{2i\sqrt{\log m}}|\CF_{\Delta_{i-1}}\right)
\\&=\p_{\left(IA_{\Delta_{i-1}}^{m,N},IA_{\Delta_{i-1}}^{m+1,N}\right)}\left(\Delta_1\le \frac{\delta N}{2i\sqrt{\log m}}\right)
\\&=\p_{\left(IA_{\Delta_{i-1}}^{m,N},IA_{\Delta_{i-1}}^{m+1,N}\right)}\left(\sum\limits_{j=1}^{ \frac{\delta N}{2i\sqrt{\log m}}}\mathbbm{1}_{\Delta_1=j}\ge 1\right),
\eae
while by \eqref{5_1} and Lemma \ref{lem2_1,1},
\bae\label{5_1,2}
&\p_{\left(IA_{\Delta_{i-1}}^{m,N},IA_{\Delta_{i-1}}^{m+1,N}\right)}\left(\sum\limits_{j=1}^{ \frac{\delta N}{2i\sqrt{\log m}}}\mathbbm{1}_{\Delta_1=j}\ge 1\right)
\\&=\p_{\left(IA_{\Delta_{i-1}}^{m,N},IA_{\Delta_{i-1}}^{m+1,N}\right)}\left(\sum\limits_{j=1}^{ \frac{\delta N}{2i\sqrt{\log m}}}\mathbbm{1}_{\Delta_1=j}\ge 1,\sum\limits_{j=1}^{ \frac{\delta N}{2i\sqrt{\log m}}}\p\left(\Delta_1=j|\CF_{j-1}\right)\le C\delta\right)
\\&\le\exp{\left[-\left(1-C\delta\right)^2/2\right]}
\\&\triangleq\delta_0
\eae
when $C\delta<1$.
It follows from \eqref{5_1,1} and \eqref{5_1,2} that
\bae\label{5_00}
\p\left(X_i=1|\CF_{\Delta_{i-1}}\right)\le\delta_0.
\eae
Again by Lemma \ref{lem2_1,1},
\bae\label{5_02}
\p\left(A_k^c\right)&=\p\left(\sum\limits_{i\in I_k}X_i\ge c_0m^{\alpha/2}\right)
\\&=\p\left(\sum\limits_{i\in I_k}X_i\ge c_0m^{\alpha/2},\sum\limits_{i\in I_k}\p\left(X_i=1|\CF_{\Delta_{i-1}}\right)\le\delta_0m^{\alpha/2}\right)
\\&\le\exp\left(-\frac{\left(c_0-\delta_0\right)^2}{c_0}m^{\alpha/2}\right).
\eae
Thus \eqref{suffice} is true by choosing adequate $c_0,\delta$, which implies \eqref{eq_lem5_1}.
\end{proof}
Now we have proved that for any $\alpha>0$, with high probability there is no more than $m^\alpha$ elements in $T^m_\Delta$. Next we want to show that all these discrepancies are highly unlikely to reach any finite subgraph $K\subseteq\vec \BH$. The proof of the following lemma is inspired by the proof of Lemma 7.1. in \cite{procaccia2019stationary}.
\begin{lemma}\label{lem5_2}
For any finite subgraph $K\subseteq\vec\BH$,
$$
\p\left(V^{D,m}_{\Delta_{m^\alpha}}\cap K\ne\emptyset\right)\le m^{-1-3\alpha/2}.
$$
\end{lemma}
\begin{proof}
For each $1\le n\le m^\alpha$, note that  
$$
\left\{\vec e_{\Delta_n,1},\vec e_{\Delta_n,2}\right\}=E^{D,m}_{\Delta_n}\backslash E^{D,m+1}_{\Delta_n-1}.
$$

For any $\vec e, A\subseteq\BZ^2$, define
\begin{align*}
&Dist\left(\vec e_1,\vec e_2\right)=\max\left\{\|\vec e_1\left(i\right)-\vec e_2\left(j\right)\|,i,j=1,2\right\},
\\&Dist\left(\vec e,A\right)=\max\left\{\|\vec e_1\left(i\right)-x\|,i=1,2,x\in A\right\}
\end{align*} with the convention that $d\left(\vec e,\emptyset\right)=\infty$.
Like  \cite{procaccia2019stationary} we have the following definitions:

\begin{itemize}
\item For any $i\ge1$, we say $\Delta_i$ is good if either $\Delta_i=\infty$ or 
$$
Dist\left(\vec e_{\Delta_i,1},\vec e_{\Delta_i,2}\right)<m^{1-5\alpha}.
$$
\item For any $i\ge1$, if $\Delta_i$ is bad, we say $\Delta_i$ is devastating if and only if $\vec e_{\Delta_i,2}$ intersects with $\left[-m^{1-3\alpha},m^{1-3\alpha}\right]\times\left[0,\log m\right]$.
 \end{itemize}
 Let
 $$
 \kappa=\inf\left\{i\ge 1:\Delta_i\text{ is bad}\right\}.
 $$
Define
\begin{itemize}
\item Event $A$: $\exists \kappa<m^\alpha$, and $\Delta_\kappa$ is devastating.\\
\item Event $B$:  $\exists \kappa<m^\alpha$, $\Delta_\kappa$ is bad but not devastating, and there is at least one bad event within $\kappa+1,\kappa+2,\cdots,m^\alpha$.
 \end{itemize}
 Then on the event $A^c\cap B^c$, for any finite $K\subseteq\vec \BL^2$,
 $$
 Dist\left(\vec e_{\Delta_i,2},K\right)\ge m^{1-3\alpha}-\sum\limits_{i=1}^{m^\alpha}m^{1-5\alpha}\ge m^{1-3\alpha}/2
 $$
 for all $1\le i\le m^\alpha$ when $m$ is large enough so that $V^{D,m}_{\Delta_{m^\alpha}}\cap K=\emptyset$.
 Thus $V^{D,m}_{\Delta_{m^\alpha}}\cap K\ne\emptyset$ implies $A$ or $B$ happens.
Define 
\begin{align*}G_k=\left\{\Delta_i \text{ is good for }i=1,\dots,k-1\right\}.
\end{align*}
Then we first present an upper bound on $\p\left(A\right)$,
\bae\label{5_03}
\p\left(A\right)&=\sum\limits_{k=1}^{m^\alpha}\p\left(G_k,\Delta_k \text{ is devastating }\right)
\\&=\sum\limits_{k=1}^{m^\alpha}\sum\limits_{j=0}^{\infty}\sum\limits_{\left(\bar A_0,\tilde A_0\right)}\p\left(G_k,\Delta_{k-1}<\infty,\Delta_k-\Delta_{k-1}>j,\left(\hat A_{\Delta_{k-1}+j}^m,\hat A_{\Delta_{k-1}+j}^{m+1}\right)=\left(\bar A_0,\tilde A_0\right)\right)
\\&\times\p_{\left(\bar A_0,\tilde A_0\right)}\left(\Delta_1=1,\Delta_1\text{ is devastating }\right)
\\&=\sum\limits_{k=1}^{m^\alpha}\sum\limits_{j=0}^{\infty}\sum\limits_{\left(\bar A_0,\tilde A_0\right)}\p\left(G_k,\Delta_{k-1}<\infty,\Delta_k-\Delta_{k-1}>j,\left(\hat A_{\Delta_{k-1}+j}^m,\hat A_{\Delta_{k-1}+j}^{m+1}\right)=\left(\bar A_0,\tilde A_0\right)\right)
\\&\times\p_{\left(\bar A_0,\tilde A_0\right)}\left(\Delta_1=1\right)\p_{\left(\bar A_0,\tilde A_0\right)}\left(\Delta_1=1,\Delta_1 \text{ is devastating }|\Delta_1=1\right).
\eae
For any $k=1,2,\dots,m^\alpha$,
\bae\label{5_55}
\p\left(G_k,\Delta_k<\infty\right)=&\sum\limits_{j=0}^{\infty}\sum\limits_{\left(\bar A_0,\tilde A_0\right)}\p\left(G_k,\Delta_{k-1}<\infty,\Delta_k-\Delta_{k-1}>j,\left(\hat A_{\Delta_{k-1}+j}^m,\hat A_{\Delta_{k-1}+j}^{m+1}\right)=\left(\bar A_0,\tilde A_0\right)\right)
\\&\times\p_{\left(\bar A_0,\tilde A_0\right)}\left(\Delta_1=1\right)\le 1,
\eae
while for any $\left(\bar A_0,\tilde A_0\right)$ satisfying 
$$
 \left\{G_k,\Delta_{k-1}<\infty,\Delta_k-\Delta_{k-1}>j,\left(\hat A_{\Delta_{k-1}+j}^m,\hat A_{\Delta_{k-1}+j}^{m+1}\right)=\left(\bar A_0,\tilde A_0\right)\right\},
 $$
 we must have $$ \bar E_0\triangle\tilde E_0\subseteq\left[\left(-\infty,-m+2m^{1-4\alpha}\right)\cup\left(m-2m^{1-4\alpha},\infty\right)\right]\times\left[0,\log m\right],$$ which is disjoint with 
 $$
 Box=\left[-2m^{1-3\alpha},2m^{1-3\alpha}\right]\times\left[0,\log m\right].
 $$
Applying Remark \ref{rem4_1} and Lemma \ref{lem2} again, we have
 \bae\label{5_3}
\p_{\left(\bar A_0,\tilde A_0\right)}\left(\Delta_1=1\right)=\sum\limits_{\vec e_{1,1}\left(2\right)\in  \bar V_0\cap\left(\tilde V_0\right)^c}\eharm_{\bar V_0\cup\tilde V_0\cup D_N}\left(\vec e_{1,1}\right)+\sum\limits_{\vec e_{1,1}\left(2\right)\in  \tilde V_0\cap\left(\bar V_0\right)^c}\eharm_{\bar V_0\cup\tilde V_0\cup D_N}\left(\vec e_{1,1}\right).
\eae

Moreover,
\bae\label{10}
&\p_{\left(\bar A_0,\tilde A_0\right)}\left(\Delta_1=1,\Delta_1 \text{ is devastating }\right)
\\&=\sum\limits_{\vec e_{1,1}\left(2\right)\in  \bar V_0\cap\left(\tilde V_0\right)^c}\eharm_{\bar V_0\cup\tilde V_0\cup D_N}\left(\vec e_{1,1}\right)\sum\limits_{{\vec e_{1,2}\left(2\right)}\in  \tilde V_0,||{\vec e_{1,2}}\left(2\right)||\le 2m^{1-3\alpha}}\eharm_{\tilde V_0\cup D_N}\left(\vec e_{1,1}\left(2\right),{\vec e_{1,2}}\right)
\\&+\sum\limits_{\vec e_{1,1}\left(2\right)\in  \tilde V_0\cap\left(\bar V_0\right)^c}\eharm_{\bar V_0\cup\tilde V_0\cup D_N}\left(\vec e_{1,1}\right)\sum\limits_{{\vec e_{1,2}\left(2\right)}\in  \bar V_0,||{\vec e_{1,2}}\left(2\right)||\le 2m^{1-3\alpha}}\eharm_{\bar V_0\cup D_N}\left(\vec e_{1,1}\left(2\right),{\vec e_{1,2}}\right)
\\&\le\sum\limits_{\vec e_{1,1}\left(2\right)\in  \bar V_0\cap\left(\tilde V_0\right)^c}\eharm_{\bar V_0\cup\tilde V_0\cup D_N}\left(\vec e_{1,1}\right)\sup_{z\in\bar V_0\Delta\tilde V_0}\sum\limits_{{\vec e_{1,2}\left(2\right)}\in  \tilde V_0,||{\vec e_{1,2}}\left(2\right)||\le 2m^{1-3\alpha}}\eharm_{\tilde V_0\cup D_N}\left(z,\vec e_{1,2}\right)
\\&+\sum\limits_{\vec e_{1,1}\left(2\right)\in  \tilde V_0\cap\left(\bar V_0\right)^c}\eharm_{\bar V_0\cup\tilde V_0\cup D_N}\left(\vec e_{1,1}\right)\sup_{z\in\bar V_0\Delta\tilde V_0}\sum\limits_{{\vec e_{1,2}\left(2\right)}\in  \bar V_0,||{\vec e_{1,2}}\left(2\right)||\le 2m^{1-3\alpha}}\eharm_{\bar V_0\cup D_N}\left( z,\vec e_{1,2}\right)
\\&\le\left(\sum\limits_{\vec e_{1,1}\left(2\right)\in  \bar V_0\cap\left(\tilde V_0\right)^c}\eharm_{\bar V_0\cup\tilde V_0\cup D_N}\left(\vec e_{1,1}\right)+\sum\limits_{\vec e_{1,1}\left(2\right)\in  \tilde V_0\cap\left(\bar V_0\right)^c}\eharm_{\bar V_0\cup\tilde V_0\cup D_N}\left(\vec e_{1,1}\right)\right)\sup_{z\in \bar V_0\Delta\tilde V_0}\p_z\left(\tau_{Box}<\tau_{D_N}\right).
\eae
Combine \eqref{5_03}, \eqref{5_55}, \eqref{5_3} and \eqref{10}, we get that
\bae\label{5_04}
\p(A)\le m^{\alpha}\sup_{z\in \bar V_0\Delta\tilde V_0}\p_z\left(\tau_{Box}<\tau_{D_N}\right).
\eae
For any $z\in V_0\Delta\tilde V_0$, since $m\le N^{1/5},$
\bae\label{11}
\p_z\left(\tau_{Box}<\tau_{D_N}\right)&= \p_z\left(\tau_{Box}<\tau_{\ell_0}\right)+\p_z\left(\tau_{\ell_0}<\tau_{Box}<\tau_{D_N}\right)
\\&\le \p_z\left(\tau_{Box}<\tau_{\ell_0}\right)+\sum\limits_{w\in D_{2N}\backslash D_N}\p_z\left(\tau_{\ell_0}=w\right)\p_w\left(\tau_{Box}<\tau_{D_N}\right)+\p_z\left(\tau_{\ell_0}<\tau_{D_{2N}}\right)
\\&\le \p_z\left(\tau_{Box}<\tau_{\ell_0}\right)+m^{1-3\alpha}\log m\sup_{w\in \ell_0\backslash D_N}\sup_{v\in Box}\p_w\left(\tau_{v}<\tau_{D_N}\right)+\frac{C\log m}{m^5},
\eae
while by Lemma 7.2 of \cite{procaccia2019stationary}, for any $\alpha<1/5$,
\bae\label{5_1,1,1}
\p_z\left(\tau_{Box}<\tau_{\ell_0}\right)\le m^{-2-3\alpha/2}.
\eae
And by the reversibility of the SRW, for any $w\in D_{2N}\backslash D_N,v\in Box,m\le N^{1/5}$, by Lemma 3.13. of \cite{procaccia2018stationary},
\bae\label{5_1,1,2}
\p_w\left(\tau_{v}<\tau_{D_N}\right)&=\frac{\p_v\left(\tau_{w}<\tau_{D_N}\right)}{\p_w\left(\tau_{w}>\tau_{D_N}\right)}
\\&\le C(\log m)^2\p_v\left(\tau_{w}<\tau_{D_N}\right)
\\&\le C(\log m)^2\p_v\left(\tau_{\partial^{out}B(v,m^5-m^{1-3\alpha})}<\tau_{D_N}\right)
\\&\le\frac{C(\log m)^3}{m^5-m^{1-3\alpha}}
\eae
for all sufficiently large $m$ since $||v-w||\ge m^5-m^{1-3\alpha}$ where $B(v,m^5-m^{1-3\alpha})$ denotes the ball centered at $v$ with radius $m^5-m^{1-3\alpha}$.  

Combine \eqref{5_04},\eqref{11}, \eqref{5_1,1,1} and \eqref{5_1,1,2}, we have 
\bae\label{12}
\p\left(A\right)\le m^\alpha\sup_{z\in\bar V_0\Delta\tilde V_0}\p_z\left(\tau_{Box}<\tau_{D_N}\right)\le m^{-1-3\alpha/2}.
\eae
Now we come to the upper bound on $\p\left(B\right)$, define
$$
B_k=\left\{\Delta_1,\cdots,\Delta_{k-1}\text{ are good},\Delta_k\text{ is bad}\right\}.
$$
Then by strong Markov property, 
$$
\p\left(B\right)\le \sum\limits_{k=1}^{m^\alpha-1}\sum\limits_{\left(\bar A_0,\tilde A_0\right)}\p\left(B_k,\Delta_k\text{ is not devastating},\left(\hat A_{\Delta_k}^m,\hat A_{\Delta_k}^{m+1}\right)=\left(\bar A_0,\tilde A_0\right)\right)
\left[\sum\limits_{j=1}^{m^\alpha-k}\p_{\left(\bar A_0,\tilde A_0\right)}\left(B_j\right)\right].
$$
Then for any configuration $\left(A,B\right)$, any $k\ge 1$,
\bae
\p_{\left(A,B\right)}\left(B_k\right)&\le \p_{\left(0,\log m\right)}\left(\tau_{\partial^{out}\left(\left[-m^{1-5\alpha}/2,m^{1-5\alpha}/2\right]\times \left[1,m^{1-5\alpha}/2\right]\right)}<\tau_{D_N}\right)
\\&\le 2\p_{\left(0,\log m\right)}\left(\tau_{\left[-m^{1-5\alpha}/2,m^{1-5\alpha}/2\right]\times \left\{m^{1-5\alpha}/2\right\}}<\tau_{D_N}\right) 
\\&\le m^{-1+6\alpha}
\eae
when $m$ is sufficiently large,
which implies that
\bae\label{13}
\p\left(B\right)\le m^{-1+7\alpha}\left[\sum\limits_{k=1}^{n^\alpha-1}\p\left(B_k\right)\right]\le m^{-2+14\alpha}.
\eae
Now by \eqref{12} and \eqref{13}, we complete the proof.
\end{proof}
%{\bf Proof of Proposition \ref{prop_coupling}:}\\
\begin{proof}[Proof of Proposition \ref{prop_coupling}]
  Combine Lemma \ref{lem5_1}, \ref{lem5_2} and Remark \ref{rem4_2}, we get Proposition \ref{prop_coupling} immediately.
%$\hfill\square$
\end{proof}

% !TEX root = Convergences_to_SDLA_3.12.19.tex

  \section{Proof of Proposition \ref{prop_coupling2}}\label{Sec_prop_coupling2}

 In this section, we consider the coupled process $\left(IA_k^{N,N^{1/5}},IA_k^{N,N}\right),k\le 2N$ whose construction is delineated in Section \ref{sec3}.
 Define
$$
\Gamma_N=\inf\left\{k\ge 1:IA_{k}^{N,N^{1/5}}\cup IA_{k}^{N,N}\not\subseteq  \vec B_1\cup \vec B_2\right\}
$$
and the truncated process before $\Gamma_N$ 
$$
\left(\hat A^{N^{1/5}}_k, \hat A^{N}_k\right)=\left(IA_{k\land \Gamma_N}^{N^{1/5},N}, IA_{k\land \Gamma_N}^{N,N}\right).
$$

By Lemma \ref{lem2_2}, there exists $C\in \left(0,\infty\right)$ such that for all sufficiently large $N$,
$$
\p\left(\Gamma_N<2N\right)<\frac{1}{N^C}.
$$
Then it suffices to show that for any $\ep>0$, there exists $N_0$ such that for any $N\ge N_0$ and any finite subgraph $K\subseteq \vec\BH$,
$$
\p\left(\exists k\le 2N,\ \hat A_k^{N^{1/5}}\cap K\ne\hat A_k^{N}\cap K\right)<\ep.
$$

   Now we divide $B_2$ into two boxes such that
  \bae
  &B_3=[-N^{1/5},N^{1/5}]\times[-\log N,\log N]\cap \BZ^2,
  \\&B_4=B_2\backslash B_3.
  \eae
  
 Since there can be too many discrepancies in $\vec B_1\cup \vec B_4$, we have to focus on the discrepancies in $\vec B_3$. Denote the vertex discrepancies set and the edge discrepancies set  constrained in $\vec B_3$ as 
 \begin{align*}
&V_n^{D,N^{1/5}}=\left\{x\in\BH: \exists k\le n \text{ s.t. } x\in \left(\hat V_k^{N^{1/5}}\Delta \hat V_k^{N}\right)\cap  B_3\right\},
 \\  &E_n^{D,N^{1/5}}=\left\{\vec e\in\vec\BH: \exists k\le n \text{ s.t. } \vec e\in \left(\hat E_k^{N^{1/5}}\Delta \hat E_k^{N}\right)\cap \left(\vec B_3\cup\widetilde{\partial^eB_3}\right)\right\}.
  \end{align*}
 
  Then we get the stopping times to creat discrepancies in $\vec B_3$ such that
  \bae
 & \Delta^{N^{1/5},N}_1=\inf\left\{1\le k\le 2N:|E_k^{D,N^{1/5}}\backslash E_{k-1}^{D,N^{1/5}}|\ge 1\right\},
 \\& \Delta^{N^{1/5},N}_i=\inf\left\{ 2N\ge k> \Delta^{m_1,m_2}_{i-1}:|E_k^{D,N^{1/5}}\backslash E_{k-1}^{D,N^{1/5}}|\ge 1\right\}
 \eae
with the convention that $\inf\emptyset=\infty$.
  Let $T_\Delta^{1/5}$ be the set of the stopping times in which a discrepancy occurs, so that
$$
T^{1/5}_\Delta=\left\{\Delta^{N^{1/5},N}_i,\Delta^{N^{1/5},N}_i\le 2N\right\}.
$$Then we want to get an upper bound on $|T^{1/5}_\Delta|$, the number of stopping times before $2N\land\Gamma_N$.

   \begin{lemma}\label{lem 6_1}
For any $\ep>0$,
$$
\p\left(|T_\Delta^{1/5}|\ge N^{2\ep}\right)\le \exp\left(-N^{\ep/2}\right).
$$for all sufficiently large $N$.
\end{lemma}
\begin{proof}
Let
$$
B_4=C_1^{\ep}\cup\cdots\cup C^{\ep}_l
$$ where $l=\lfloor\frac{N^{1-\ep/2}}{2}-N^{1/5-\ep/2}\rfloor$, $C^\ep_i=[N^{1/5}+\left(i-1\right)N^{\ep/2},N^{1/5}+iN^{\ep/2}]\times [-\log N,\log N]\cap \BZ^2,1\le i\le l-1$, and $C^\ep_l=[N^{1/5}+\left(l-1\right)N^{\ep/2},N/2]\times [-\log N,\log N]\cap \BZ^2$.
Thus we divide $B_1\cup B_2$ into $l+2$ parts such that 
\bae\label{6_00}
 \hat V_n^{N^{1/5}}\cup\hat V_n^{N}\subseteq B_1\cup B_3\cup \left(\cup_{1\le i\le l}C^\ep_i\right).
 \eae
 By Remark \ref{rem4_1}, Lemma \ref{lem2}, \eqref{6_00},
\bae\label{6_66}
 & \p\left(n+1\in T_\Delta^{1/5}|\CF_{n}\right)
 \\&=\sum\limits_{\vec e_{n,1}\left(2\right)\in  \hat V_n^{N^{1/5}}\cap\left(\hat V_n^{N}\right)^c\cap \left(B_1\cup B_4\right)}\eharm_{\hat V_n^{N^{1/5}}\cup \hat V_n^{N}\cup D_N}\left(\vec e_{n,1}\right)\sum\limits_{{\vec e_{n,2}\left(2\right)}\in \hat V_n^{N}\cap B_3}\eharm_{\hat V_n^{N}\cup D_N}\left(\vec e_{n,1}\left(2\right),{\vec e_{n,2}}\right)
\\&+\sum\limits_{\vec e_{n,1}\left(2\right)\in  \hat V_n^{N}\cap\left(\hat V_n^{N^{1/5}}\right)^c\cap \left(B_1\cup B_4\right)}\eharm_{\hat V_n^{N^{1/5}}\cup \hat V_n^{N}\cup D_N}\left(\vec e_{n,1}\right)\sum\limits_{{\vec e_{n,2}\left(2\right)}\in \hat V_n^{N^{1/5}}\cap B_3}\eharm_{\hat V_n^{N}\cup D_N}\left(\vec e_{n,1}\left(2\right),{\vec e_{n,2}}\right)
\\&+\sum\limits_{\vec e_{n,1}\left(2\right)\in \left( \hat V_n^{N^{1/5}}\Delta\hat V_n^{N}\right)\cap B_3}\eharm_{\hat V_n^{N^{1/5}}\cup \hat V_n^{N}\cup D_N}\left(\vec e_{n,1}\right)
  \\&\le \sum\limits_{\vec e\left(2\right)\in B_1}\eharm_{B_1\cup D_N}\left(\vec e\right)\sup_{z\in B_1}\p_z\left(\tau_{B_3}<\tau_{D_N}\right)
  \\&+\sum\limits_{n=1}^{\lfloor\frac{N^{1-\ep/2}}{2}-N^{1/5-\ep/2}\rfloor}\sum\limits_{\vec e\left(2\right)\in C_n^\ep}\eharm_{C_n^\ep\cup D_N}\left(\vec e\right)\sup_{z\in C_n^\ep}\p_z\left(\tau_{B_3}<\tau_{D_N}\right)
 +\frac{C\sqrt{\log N}}{N}|V_{n}^{D,N^{1/5}}|
  \\&=I_1+I_2+I_3.
 \eae

For $I_1$,
\bae\label{6_1}
I_1&=\sum\limits_{\vec e\in \partial^eB_1}\eharm_{B_1\cup D_N}\left(\vec e\right)\sup_{z\in B_1}\p_z\left(\tau_{B_3}<\tau_{D_N}\right)
\\ &\le\sup_{z\in B_1}\p_z\left(\tau_{B_3}<\tau_{D_N}\right)
 \\&\le\sup_{z\in B_1}[\p_z\left(\tau_{H_{N/4}}<\tau_{D_N},||S_{\tau_{H_{N/4}}}||\ge N^4\right)
 \\&+\sum\limits_{||w||\le N^4}\p_z\left(\tau_{H_{N/4}}<\tau_{D_N},S_{\tau_{H_{N/4}}}=w\right)\p_w\left(\tau_{B_3}<\tau_{D_N}\right).
 \eae
And by the Beurling estimate, Theorem 1 of \cite{lawler2004beurling},
 \bae\label{6_8}
 &\sup_{z\in B_1}\p_z\left(\tau_{H_{N/4}}<\tau_{D_N},||S_{\tau_{H_{N/4}}}||\ge N^4\right)
  \\&\le\p_0\left(\tau_{N^4/2}<\tau_{A[\sqrt{N},N^4/4]}\right)
 \\&\le c\sqrt{\frac{2}{N^{4-1/2}}},
\eae  where the second inequality comes from 
\begin{align*}
&\left\{y\in H_{N/4},||y||\ge N^4\right\}\subseteq \partial^{out}B\left(z,N^4/2\right),
\\&B\left(z,N^4/4\right)\backslash B\left(z,\sqrt{N}\right)\subseteq D_N\cup \left\{z\in H_{N/4},||y||\le N^4\right\}
\end{align*}for each $z\in B_1$.
Moreover, 
\bae\label{6_9}
 &\sup_{z\in B_1}\sum\limits_{||w||\le N^4}\p_z\left(\tau_{H_{N/4}}<\tau_{D_N},S_{\tau_{H_{N/4}}}=w\right)\p_w\left(\tau_{B_3}<\tau_{D_N}\right)
 \\&\le \sup_{z\in B_1} \p_z\left(\tau_{H_{N/4}}<\tau_{D_N}\right)\sup_{||w||\le N^4,w\in H_{N/4}}\p_w\left(\tau_{B_3}<\tau_{D_N}\right)
 \\&\le N^{1/5}\log N\sup_{z\in B_1}\p_z\left(\tau_{H_{N/4}}<\tau_{D_N}\right)\sup_{||w||\le N^4,w\in H_{N/4}}\sup_{\tilde w\in B_3}\p_w\left(\tau_{\tilde w}<\tau_{D_N}\right)
 \\&\le\frac{\log N}{N^{1/20}}\sup_{||w||\le N^4,w\in H_{N/4}}\sup_{\tilde w\in B_3}\p_w\left(\tau_{\tilde w}<\tau_{D_N}\right),
 \eae
 while for any $||w||\le N^4,w\in H_{N/4},\tilde w\in B_3$, by the reversibility of the SRW, Lemma 3.13. of \cite{procaccia2018stationary}, and the Beurling estimate,
 \bae\label{6_01}
 \p_w\left(\tau_{\tilde w}<\tau_{D_N}\right)&=\frac{\p_{\tilde w}\left(\tau_{ w}<\tau_{D_N}\right)}{\p_w\left(\tau_{w}>\tau_{D_N}\right)}
\le \frac{C(\log N)^3}{N}
\eae
 
 Combine \eqref{6_1}, \eqref{6_8}, \eqref{6_9} and \eqref{6_01}, we have
 \bae\label{6_2}
 I_1=o\left( \frac{1}{N}\right).
 \eae

 For $I_2$, similarly, by Lemma \ref{lem2},
 \bae\label{6_3}
I_2=&\sum\limits_{n=1}^{\lfloor\frac{N^{1-\ep/2}}{2}-N^{1/5-\ep}\rfloor}\sum\limits_{\vec e\in \partial^eC_n^{\ep}}\eharm_{C_n^{\ep}\cup D_N}\left(\vec e\right)\sup_{z\in C_n^{\ep}}\p_z\left(\tau_{B_3}<\tau_{D_N}\right)
\\&\le\sum\limits_{n=1}^{\lfloor\frac{N^{1-\ep/2}}{2}-N^{1/5-\ep/2}\rfloor}N^\ep \log N \frac{C\sqrt{\log N}}{N} \sup_{z\in C_n^{\ep}}\p_z\left(\tau_{B_3}<\tau_{D_N}\right)
\\&\le \frac{C\left(\log N\right)^2}{N^{1-\ep/2}}\left(\sum\limits_{n=1}^{\lfloor\frac{N^{1-\ep/2}}{2}-N^{1/5-\ep/2}\rfloor} \sup_{z\in C_n^{\ep}}\p_z\left(\tau_{B_3}<\tau_{D_N}\right)\right)
\\&\le \frac{C\left(\log N\right)^2}{N^{1-\ep/2}}\left(\sum\limits_{n=1}^{N^{1-\ep/2}}\frac{\log N}{nN^{\ep/2}}+1\right)
 \\& \le\frac{C\left(\log N\right)^3}{N^{1-\ep/2}}.
\eae 
 Thus it follows from \eqref{6_66}, \eqref{6_2} and \eqref{6_3} that
 \bae\label{6_04}
  \p\left(n\in T_\Delta^{1/5}|\CF_{n-1}\right)\le \frac{1}{N^{1-\ep}}+\frac{C\sqrt{\log N}}{N}|V_{n-1}^{D,N^{1/5}}|.
 \eae
 Applying the proof  in Lemma \ref{lem5_1}, for any $\delta<1$, let $\Delta_0=0$ and
$$
\forall 1\le i\le N^{2\ep},X_i=
\begin{cases}
1&\text{ if $\Delta_i-\Delta_{i-1}\le \frac{\delta N}{2i\sqrt{\log N}+N^\ep}$ or $\Delta_i=\infty$}\\
0&\text{ otherwise}
\end{cases}.
$$
Define
\begin{align*}
&\forall 1\le k\le N^\ep,I_k=\left\{\left(k-1\right)N^\ep+1,\cdots,kN^\ep\right\},
\\&A_k=\left\{\sum\limits_{i\in I_k}X_i<c_0N^\ep\right\}.
\end{align*}
On $\cap_{1\le k\le N^\ep}A_k$,
$$
\sum\limits_{i=1}^{N^{2\ep}}\left(\Delta_i-\Delta_{i-1}\right)\ge\sum\limits_{k=1}^{N^{\ep}}\left(c_0N^\ep \times\frac{\delta N}{2kN^\ep\sqrt{\log N}+N^\ep}\right)\ge\frac{\ep c_0\delta N\sqrt{\log N}}{4}\ge 2N
$$for any $c_0,\delta>0$ when $N$ is sufficiently large.
So that it suffices to show that for any $1\le k\le N^{2\ep}$,
$$
\p\left(A_k^c\right)\le  \exp\left(-CN^{\ep}\right).
$$
Notice that
\bae\label{6_1,1,1}
\p\left(X_i=1|\CF_{\Delta_{i-1}}\right)&=\p\left(\Delta_i-\Delta_{i-1}\le \frac{\delta N}{2i\sqrt{\log N}+N^\ep}|\CF_{\Delta_{i-1}}\right)
\\&=\p_{\left(IA_{\Delta_{i-1}}^{m,N},IA_{\Delta_{i-1}}^{m+1,N}\right)}\left(\Delta_1\le \frac{\delta N}{2i\sqrt{\log N}+N^\ep}\right)
\\&=\p_{\left(IA_{\Delta_{i-1}}^{m,N},IA_{\Delta_{i-1}}^{m+1,N}\right)}\left(\sum\limits_{j=1}^{ \frac{\delta N}{2i\sqrt{\log N}+N^\ep}}\mathbbm{1}_{\Delta_1=j}\ge 1\right),
\eae
while by \eqref{6_04} and Lemma \ref{lem2_1,1},
\bae\label{6_1,2,1}
&\p_{\left(IA_{\Delta_{i-1}}^{m,N},IA_{\Delta_{i-1}}^{m+1,N}\right)}\left(\sum\limits_{j=1}^{ \frac{\delta N}{2i\sqrt{\log N}+N^\ep}}\mathbbm{1}_{\Delta_1=j}\ge 1\right)
\\&=\p_{\left(IA_{\Delta_{i-1}}^{m,N},IA_{\Delta_{i-1}}^{m+1,N}\right)}\left(\sum\limits_{j=1}^{ \frac{\delta N}{2i\sqrt{\log N}+N^\ep}}\mathbbm{1}_{\Delta_1=j}\ge 1,\sum\limits_{j=1}^{ \frac{\delta N}{2i\sqrt{\log N}+N^\ep}}\p\left(\Delta_1=j|\CF_{j-1}\right)\le C\delta\right)
\\&\le\exp{[-\left(1-C\delta\right)^2/2]}
\\&\triangleq\delta_0
\eae
when $C\delta<1$.
It follows from \eqref{6_1,1,1} and \eqref{6_1,2,1} that
\bae\label{6_05}
\p\left(X_i=1|\CF_{\Delta_{i-1}}\right)\le\delta_0.
\eae
Again by \eqref{6_05},
\bae
\p\left(A_k^c\right)&=\p\left(\sum\limits_{i\in I_k}X_i\ge c_0N^\ep\right)
\\&=\p\left(\sum\limits_{i\in I_k}X_i\ge c_0N^\ep,\sum\limits_{i\in I_k}\p\left(X_i=1|\CF_{\Delta_{i-1}}\right)\le\delta_0N^\ep\right)
\\&\le\exp\left(-\frac{\left(c_0-\delta_0\right)^2}{c_0}N^\ep\right).
\eae
Thus by choosing adequate $c_0,\delta$, we have
 $$
\p\left(|T_\Delta^{1/5}|\ge N^{2\ep}\right)\le \exp\left({-N^{\ep/2}}\right).
$$

 \end{proof}

 \begin{proposition}
\label{prop5_1}
For any finite subgraph $K\subseteq \vec\BL^2,$ any $\ep>0$, there exists $N_0>0$ such that for all $N\ge N_0$, 
	\bae
			\prob\left(\hat A^{N^{1/5}}_k\cap K=  \hat A^{N}_k\cap K, \ \forall k\le 2N\right)\ge 1-\ep.
	\eae
\end{proposition}
\begin{proof} 
For any $i\ge1$, we say $\Delta_i$ is good if either $\Delta_i=\infty$ or 
$$
Dist\left(\vec e_{\Delta_i,1},\vec e_{\Delta_i,2}\right)<N^{1/10-2\ep}.
$$
Define
\begin{itemize}
\item Event $A$:  $\exists\ \Delta_i\in T_\Delta^{1/5}$ such that $\Delta_i$ is bad. 
\item Event $B$: $\exists \ \Delta_i\in T_\Delta^{1/5}$ such that $\vec e_{\Delta_i,1}\left(2\right)\in B_1\cup B_4$ and $\vec e_{\Delta_i,2}\left(2\right)\in B_5$.\end{itemize}
It is easy to see that
\bae\label{6_17}
&A^c\cap B^c\cap \left\{|T_\Delta^{1/5}|\le N^\ep\right\}
\\&\subseteq \left\{E^{D,N^{1/5}}_{2N}\subseteq \left( [N^{1/10}-N^{1/10-\ep},N^{1/5}]\cup[-N^{1/5},-N^{1/10}+N^{1/10-\ep}]\right)\times [-\log N,\log N]\right\}
\\&\subseteq \left\{V^{D,N^{1/5}}_{2N}\cap K=\emptyset\right\}.
\eae
Thus by Lemma \ref{lem 6_1} and \eqref{6_17}, we have
\bae\label{6_8,1}
& \prob\left(\exists k\le 2N,IA^{N^{1/5},N}_k\cap K\ne  IA^{N,N}_k\cap K\right)
\\&=\p\left(V_{2N}^{D,N^{1/5}}\cap K\ne\emptyset\right)
 \\&\le\p\left(|T_\Delta^{1/5}|\ge N^\ep\right)+\prob\left(V_{2N}^{D,N^{1/5}}\cap K\ne\emptyset,|T_\Delta^{1/5}|\le N^\ep\right)
 \\&\le \exp\left({-N^{\ep/2}}\right)+\p\left(|T_\Delta^{1/5}|\le N^\ep,B\right)
 +\p\left(|T_\Delta^{1/5}|\le N^\ep,A\right).
 \eae
When restricted on $B$,
 \bae\label{6_6}
\p\left(|T_\Delta^{1/5}|\le N^\ep,B\right)&\le \p\left(\exists n\le 2N,\vec e_{n,1}\left(2\right)\in B_1\cup B_4,\vec e_{n,2}\left(2\right)\in B_5\right)
\\&\le 2N\bigg[\sum\limits_{\vec e\left(2\right)\in B_1}\eharm_{B_1\cup D_N}\left(\vec e\right)\sup_{z\in B_1}\p_z\left(\tau_{B_5}<\tau_{D_N}\right)
  \\&+\sum\limits_{n=1}^{\lfloor\frac{N^{1-\ep/2}}{2}-N^{1/5-\ep/2}\rfloor}\sum\limits_{\vec e\left(2\right)\in C_n^\ep}\eharm_{C_n^\ep\cup D_N}\left(\vec e\right)\sup_{z\in C_n^\ep}\p_z\left(\tau_{B_5}<\tau_{D_N}\right)\bigg]
  \\&=2N\left(I'_1+I'_2\right).
 \eae
 For $I'_1$, since $B_5\subseteq B_3$
 \bae\label{6_1,3}
 I_1'\le I_1= o\left(\frac{1}{N}\right).
 \eae
 For $I'_2$, by Lemma \ref{lem2}.
 \bae\label{6_1,1}
\sum\limits_{\vec e\left(2\right)\in C_n^\ep}\eharm_{C_n^\ep\cup D_N}\left(\vec e\right)\le CN^{\ep/2}  \frac{\sqrt{\log N}}{N}
\eae 
And for any $1\le n\le \lfloor\frac{N^{1-\ep/2}}{2}-N^{1/5-\ep/2}\rfloor$, just as before, we have
\bae\label{6_1,2}
&\sup_{z\in C_n^{\ep}}\p_z\left(\tau_{B_5}<\tau_{D_N}\right)
 \\&\le\sup_{z\in C_n^{\ep}}\bigg[\p_z\left(\tau_{H_{N^{1/5}}}<\tau_{D_N},||S_{\tau_{H_{N^{1/5}}}}||\ge N^4\right)
 \\&+\sum\limits_{||w||\le N^4}\p_z\left(\tau_{H_{N^{1/5}}}<\tau_{D_N},S_{\tau_{H_{N^{1/5}}}}=w\right)\p_w\left(\tau_{B_5}<\tau_{D_N}\right)\bigg]
  \\&\le o\left(\frac{1}{N^{3/2}}\right)+\sup_{z\in C_n^{\ep}}\sum\limits_{||w||\le N^4}\p_z\left(\tau_{H_{N^{1/5}}}<\tau_{D_N},S_{\tau_{H_{N^{1/5}}}}=w\right)\p_w\left(\tau_{B_5}<\tau_{D_N}\right) 
   \\&\le o\left(\frac{1}{N^{3/2}}\right)+ \sup_{z\in C_n^{\ep}} \p_z\left(\tau_{H_{N^{1/5}}}<\tau_{D_N}\right)\sup_{||w||\le N^4,w\in H_{N^{1/5}}}\p_w\left(\tau_{B_5}<\tau_{D_N}\right)
 \\&\le o\left(\frac{1}{N^{3/2}}\right)+ N^{1/10}\log N\sup_{z\in C_n^{\ep}}\p_z\left(\tau_{H_{N^{1/5}}}<\tau_{D_N}\right)\sup_{||w||\le N^4,w\in H_{N^{1/5}}}\sup_{\tilde w\in B_5}\p_w\left(\tau_{\tilde w}<\tau_{D_N}\right)
 \\&\le o\left(\frac{1}{N^{3/2}}\right)+\frac{C\left(\log N\right)^3}{nN^{\ep/2+1/10}}.
\eae
It follows from \eqref{6_1,1} and \eqref{6_1,2} that
\bae\label{6_1,4}
I_2'&=\sum\limits_{n=1}^{\lfloor\frac{N^{1-\ep/2}}{2}-N^{1/5-\ep/2}\rfloor}\sum\limits_{\vec e\left(2\right)\in C_n^\ep}\eharm_{C_n^\ep\cup D_N}\left(\vec e\right)\sup_{z\in C_n^\ep}\p_z\left(\tau_{B_5}<\tau_{D_N}\right)
\\&\le N^{\ep/2} \log N \frac{C\sqrt{\log N}}{N}\times\sum\limits_{n=1}^{\lfloor\frac{N^{1-\ep/2}}{2}-N^{1/5-\ep/2}\rfloor}\left[ o\left(\frac{1}{N^{3/2}}\right)+\frac{C\left(\log N\right)^3}{nN^{\ep/2+1/10}}\right]
\\&= o\left(\frac{1}{N}\right).
\eae
When restricted on $A$,
 \bae\label{6_7}
 \p\left(|T_\Delta^{1/5}|\le N^\ep,A\right)=\p\left(|T_\Delta^{1/5}|\le N^\ep,\exists \Delta_n\in T_\Delta^{1/5},\Delta_n\text{ is bad}\right)\le N^\ep \times\frac{C\log N}{N^{1/10-2\ep}}.
 \eae
 Substitute \eqref{6_6}, \eqref{6_1,3}, \eqref{6_1,4} and \eqref{6_7} into \eqref{6_8,1}, we can get that for all sufficiently large $N$,
 $$
  \prob\left(\exists k\le N,IA^{N^{1/5},N}_k\cap K\ne  IA^{N,N}_k\cap K\right)< \ep.
 $$
 
 \end{proof}
 {\bf Proof of Proposition \ref{prop_coupling2}:}\\
 Proposition \ref{prop_coupling2} follows from Proposition \ref{prop5_1} and Remark \ref{rem4_2}.
 $\hfill\square$

% !TEX root =Convergences_to_SDLA_3.12.19.tex

\section{Appendix}\label{appendix}
\subsection{Proof of Lemma \ref{lem1_1}}
\begin{proof}
Obviously, the weak convergence implies the finite dimensional distribution's convergence. So we only need to prove the other direction.
For convenience, let
$$
X_n\left(t\right)\triangleq EA_{nt}^n\cap\vec\BH,X_\infty\left(t\right)\triangleq SA_{ct}^\infty.
$$
 Since $\left(E,\rho\right)$ is a complete and totally bounded metric space, which implies that it is also separable and compact. So that the set of the probability measures on $E$ is compact. By Theorem 7.8, (b) of \cite{ethier2009markov}, \eqref{1_1} implies the convergence of the finite dimensional distribution. In order to prove the weak convergence, by Theorem 7.8, (b) of \cite{ethier2009markov} again, we only need to prove that $\{X_n\left(t\right)\}^\infty_{n=1}$ is relatively compact. I.e. each sequence of $\{X_n\left(t\right)\}^\infty_{n=1}$ has a weakly convergent subsequence. 

Define
$$
w'\left(X_n,\delta,T\right)=\inf_{\{t_i\}}\max_{0\le i< r}\sup_{s,t\in [t_i,t_{i+1})}\rho\left(X_n\left(s\right),X_n\left(t\right)\right) 
,$$

 where $\{t_i\}$ ranges over all partitions of the form $0=t_0<t_1<\cdots<t_{r-1}<T\le t_r$ with $t_i-t_{i-1}>\delta$ for all $1\le i\le r$.
Then by Corollary 7.4 of \cite{ethier2009markov}, a necessary and sufficient condition for the relative compactness of $\{X_n\left(t\right)\}^\infty_{n=1}$ is that for each $\eta>0$ and $T\in\left(0,\infty\right)$, there exists $\delta>0$ such that
\bae\label{8_4}
\limsup_{n\to\infty}\p\left(w'\left(X_n,\delta,T\right)>\eta\right)<\eta.
\eae 

Recall the definition of $\rho$ in Section 4.1. of \cite{liggett2010continuous} such that for any $\eta,\zeta\in E$,
$$
\rho(\eta,\zeta)=\sum\limits_{x\in\vec\BH,\ x\text{ is an edge or a vertex }}\alpha(x)|\eta(x)-\zeta(x)|.
$$ 

Since $\alpha\left(x\right)$ is summable, for any $\eta>0$, there exists a finite subgraph $F\subseteq\vec\BH$ such that
$$
\sup_{\xi\equiv\zeta \text{ on } F}\rho\left(\xi,\zeta\right)\le\sum\limits_{x\in \vec\BH\backslash F , \ x\text{ is an edge or a vertex }}\alpha\left(x\right)<\eta/3,
$$
and denote
\bae\label{8_5}
M_F=\sup_{x\in F,\ x\text{ is an edge or a vertex }}\alpha\left(x\right).
\eae
For any configuration $\xi\in E$, let
$$\xi^F\left(x\right)=
\begin{cases}
\xi\left(x\right)&x\in F, x\text{ is an edge or a vertex }\\
0&\text{otherwise}.
\end{cases}
$$

Then for any $n$, by the triangle inequality of $\rho$, \eqref{8_5}, and the increasing property of $X_n^F\left(t\right)$ with respect to $t$,
\bae\label{8_1}
&\p\left(\inf_{\{t_i\}}\max_{0\le i< r}\sup_{s,t\in [t_i,t_{i+1})}\rho\left(X_n\left(s\right),X_n\left(t\right)\right) >\eta\right)
\\&\le\p\left(\inf_{\{t_i\}}\max_{0\le i< r}\sup_{s,t\in [t_i,t_{i+1})}\rho\left(X^F_n\left(s\right),X^F_n\left(t\right)\right) >\eta/3\right)
\\&\le\p\left(\inf_{\{t_i\}}\max_{0\le i< r}\sup_{s,t\in [t_i,t_{i+1})}|X_n^F\left(t\right)-X_n^F\left(s\right)|>\frac{\eta}{3M_F}\right)
\\&\le\p\left(\inf_{\{t_i\}}\max_{0\le i< r}|X_n^F\left(t_{i+1}-\right)-X_n^F\left(t_{i}\right)|>0\right)
\eae

where $$|\eta-\zeta|=\sum\limits_{x\in\vec\BH,\ x\text{ is an edge or a vertex }}|\eta\left(x\right)-\zeta\left(x\right)|.$$
Define stopping times
$$
\tau_0^n=0, \tau_k^n=\inf\{T\ge t>\tau_{k-1}^n,|X_n^F\left(t\right)-X_n^F\left(t-\right)|\ge 1\},k\ge 1
$$
with the convention that $\inf\emptyset=\infty$.

 Then on the event $\{\inf_{\{t_i\}}\max_{0\le i< r}|X_n^F\left(t_{i+1}-\right)-X_n^F\left(t_{i}\right)|>0\}$, there must be a waiting time $\Delta_k^n= \tau_k^n- \tau_{k-1}^n$ smaller than $2\delta$. Otherwise by choosing $\{t_i=\tau_i^n,i< r= N_n\left(T\right), t_{r}=T\}$, we can get a contradiction since
$$t_i-t_{i-1}>\delta, \text{ and }
\max_{0\le i< r}|X_n^F\left(t_{i+1}-\right)-X_n^F\left(t_{i}\right)|=0.
$$
So that by \eqref{8_1},
\bae\label{8_8}
\p\left(\inf_{\{t_i\}}\max_{0\le i< r}|X_n^F\left(t_{i+1}-\right)-X_n^F\left(t_{i}\right)|>\eta\right)
\le\p\left(\exists \Delta_k^n\text{ s.t. }\Delta_k^n<2\delta\right)
\eae
By Lemma \ref{lem2}, there exists a constant $C_F\in \left(0,\infty\right)$, only depend on $F$ such that 
$$
\sum\limits_{x\in F}n\eharm_{X_n\left(t\right)}\left(x\right)\le C_F
$$
 for all $t\le T$ and sufficient large $n$.
Therefore, for each $n$, $\{\tau_0^n=0, \tau_k^n\le T\}$ can be stochastically dominated by a Poisson flow $\{\tau_0^F=0, \tau_k^F\le T\}$ with intensity $C_F$. Denote the waiting times as $ \Delta_k^F=\tau_k^F-\tau_{k-1}^F$ and the number of arrivals before time $T$ as $N^F\left(T\right)$. Since conditional on $\{N^F\left(T\right)=k\}$, each arrival time is uniformly distributed on $[0,T]$, 
\bae\label{8_2}
&\p\left(\exists \Delta_k^n\text{ s.t. }\Delta_k^n<2\delta\right)
\\&\le\p\left(\exists \Delta_k^F,\text{ s.t. }\Delta_k^F<2\delta\right)
\\&= \sum\limits_{k=2}^\infty k\left(k-1\right)\p\left(N^F\left(T\right)=k\right)\p\left(0<\tau_2^F-\tau_{1}^F<2\delta|N^F\left(T\right)=k\right)
\\&\le 2C^2_F\delta T.
\eae
Then for each $\eta$, we can choose $\delta=\frac{\eta}{2C^2_FT}$ so that \eqref{8_4} comes from \eqref{8_1}, \eqref{8_8} and \eqref{8_2}.

\end{proof}

\subsection{Proof of Lemma \ref{lem2}}
\begin{proof}
Recalling the definition of the edge harmonic measure, for any $x\in\partial^{out}A$,
$$
\eharm_{A\cup D_N}\left(x\right)=\sum_{\vec e: \ \vec e\left(1\right)=x} \eharm_{A\cup D_N}\left(\vec e\right)\le \sum_{\vec e: \ \vec e\left(1\right)=x} \harm_{A\cup D_N}\left(\vec e\left(2\right)\right).
$$

Then it suffices to show that for any $\vec e\left(2\right)=y$ where $y$ is a neighbor of $ x$, 
$$
N\harm_{A\cup D_N}\left(y\right)\le C\sqrt{|y{\left(2\right)}|+1}.
$$
Without loss of generality, we can assume that $y{\left(2\right)}=n$. Since $A$ is connected and $A\cap l_0\ne\emptyset$, there must be a finite nearest neighbor path 
$$\CP_n=\{y=P_0,P_1,\cdots,P_{n_y}\in l_0\},||P_i-P_{i+1}||=1,0\le i\le n_y$$
 from $y$ to $l_0$. Since $y\left(2\right)=n$, we have $||y-P_{n_y}||\ge n$.

 Define 
 \begin{align*}
 &m_n=\inf\{i:||P_i-x||\ge n\},
 \\&Q_n=\{P_0,\cdots,P_{m_n}\},
\\&\hat P_n=Q_n\cup D_N.
\end{align*}

Then
\bae\label{5}
\harm_{A\cup D_N}\left(y\right)&\le \harm_{D_N\cup \hat P_n }\left(y\right)
\\&=\lim\limits_{R\to\infty}\frac{1}{|\partial^{out}B\left(0,R\right)|}\sum\limits_{z\in\partial^{out}B\left(0,R\right)}\harm_{D_N\cup \hat P_n }\left(z,y\right)
\\&=\lim\limits_{R\to\infty}\frac{1}{|\partial^{out}B(0,R)|}\e_y\left[\sharp \text{ visits to }\partial^{out}B(0,R) \text{ in } [0,\tau_{D_N\cup \hat P_n })\right]
\\&\le\lim\limits_{R\to\infty}\frac{C}{R}\e_y\left[\sharp \text{ visits to }\partial^{out}B(0,R) \text{ in } [0,\tau_{D_N\cup \hat P_n })\right].
\eae
Next we want to show that 
$$
\e_y[\sharp \text{ visits to } \partial^{out}B\left(0,R\right) \text{ in }[0,\tau_{D_N\cup \hat P_n })]\le C R \p_y\left(\tau_{2N}<\tau_{D_N\cup \hat P_n }\right).
$$
Since $C_N=[-\lfloor N/2\rfloor,0]\times 0\subseteq D_N$,
\bae\label{6}
&\e_y[\sharp \text{ visits to } \partial^{out}B\left(0,R\right) \text{ in }[0,\tau_{D_N\cup \hat P_n })]
\\&\le\frac{\p_y\left(\tau_R<\tau_{D_N\cup \hat P_n }\right)}{\min_{z\in\partial^{out}B\left(0,R\right)}\p_z\left(\tau_R>\tau_{D_N\cup \hat P_n }\right)}
\\&=\frac{1}{\min_{z\in\partial^{out}B\left(0,R\right)}\p_z\left(\tau_R>\tau_{D_N\cup \hat P_n }\right)}
\left[\sum\limits_{z\in\partial^{out}B\left(0,2N\right)}\p_y\left(\tau_{2N}<\tau_{D_N\cup \hat P_n },S_{\tau_{2N}}=z\right)\p_z\left(\tau_R<\tau_{D_N\cup \hat P_n }\right)\right]
\\&\le\frac{1}{\min_{z\in\partial^{out}B\left(0,2N\right)}\p_z\left(\tau_R>\tau_{D_N\cup \hat P_n }\right)}
\left[\sum\limits_{z\in\partial^{out}B\left(0,R\right)}\p_y\left(\tau_{2N}<\tau_{D_N\cup \hat P_n },S_{\tau_{2N}}=z\right)\p_z\left(\tau_R<\tau_{C_N }\right)\right]
\\&\le\frac{\p_y\left(\tau_{2N}<\tau_{D_N\cup \hat P_n }\right)\max_{z\in\partial^{out}B\left(0,2N\right)}\p_z\left(\tau_R<\tau_{C_N}\right)}{\min_{z\in\partial^{out}B\left(0,R\right)}\p_z\left(\tau_R>\tau_{D_N\cup \hat P_n }\right)}.
\eae
  
While by Lemma 3-4 of \cite{kesten1987hitting}, if $D_N\cup \hat P_n\subseteq B\left(0,r\right)$ for some $2r+1<R$,
\bae\label{3}
\min_{z\in\partial^{out}B\left(0,R\right)}\p_z\left(\tau_R>\tau_{D_N\cup \hat P_n }\right)\ge C\left(R\log R\right)^{-1},
\eae
and
\bae\label{4}
\max_{z\in\partial^{out}B\left(0,2N\right)}\p_z\left(\tau_R<\tau_{C_N}\right)\le C\left(\log R\right)^{-1}.
\eae
It follows from \eqref{5}, \eqref{6}, \eqref{3} and \eqref{4} that
\bae\label{7}
\CH_{A\cup D_N}\left(y\right)\le C\p_y\left(\tau_{2N}<\tau_{D_N\cup \hat P_n }\right).
\eae

Then we only need to show that 
\bae\label{1}
\p_y\left(\tau_{2N}<\tau_{D_N\cup \hat P_n }\right)\le \frac{Cn^{1/2}}{N}.
\eae
Define $r_n=2n,n\le \log m, S_n=\partial^{out}B\left(y,Cr_n\right)\cap\{\left(x,y\right)\in\BZ^2,y\ge 1\}\subseteq B\left(0,2N\right)$ for some proper constant $C$, so that 
\bae\label{15}
\p_y\left(\tau_{2N}<\tau_{D_N\cup \hat P_n }\right)=\sum\limits_{z\in S_n}\p_y\left(\tau_{S_n}<\tau_{D_N\cup \hat P_n },S_{\tau_{S_n}}=z\right)\p_z\left(\tau_{2N}<\tau_{D_N\cup \hat P_n }\right),
\eae
On one hand, for any $z\in S_n$, $|z\left(1\right)|\le m+\log m+2r_n$, so that when $N$ is large enough, $[z\left(1\right)-\delta N/2,z\left(1\right)+\delta N/2]\times [0,\delta N/2]\subseteq B\left(0,N\right),$ which implies 
\bae\label{16}
\p_z\left(\tau_{2N}<\tau_{D_N\cup \hat P_n }\right)\le C\p_z\left(\tau_{[z{\left(1\right)}-\delta N/2,z{\left(1\right)}+\delta N/2]\times \{\delta N/2\}}<\tau_{D_N }\right)\le Cn/N.
\eae
On the other hand, by $\left(52\right)$ of \cite{procaccia2019stationary}, 
\bae\label{17}
\p_y\left(\tau_{S_n}<\tau_{D_N\cup \hat P_n }\right)\le Cn^{-1/2}.
\eae
Now  \eqref{1} can be derived from \eqref{15}, \eqref{16} and \eqref{17}.
\end{proof}

\subsection{Proof of Lemma \ref{lem2-1}}
\begin{proof}
We will prove the result by induction. First when $n=1$, for any increasing function $f$ on $\{0,1\}$,
\bae
\e f\left(X_1\right)&=f\left(0\right)\p\left(X_1=0\right)+f\left(1\right)\p\left(X_1=1\right)
\\&=f\left(0\right)+\left[f\left(1\right)-f\left(0\right)\right]\p\left(X_1=1\right)
\\&\le f\left(0\right)+[f\left(1\right)-f\left(0\right)]\p\left(Y_1=1\right)\\&=\e f\left(Y_1\right).
\eae
Now we assume the result is true for all $n\le N-1$. We come to the case $n=N$. For any increasing function $f$ on $\{0,1\}^N$, any $\left(a_1,\cdots,a_{N}\right)\in\{0,1\}^{N}$,

\bae
&\e f\left(X_1,\cdots,X_N\right)\\&=\sum\limits_{a_1,\cdots,a_{N-1}}\p\left(X_1=a_1,\cdots,X_{N-1}=a_{N-1},X_N=0\right)f\left(a_1,\cdots,a_{N-1},0\right)
\\&+\sum\limits_{a_1,\cdots,a_{N-1}}\p\left(X_1=a_1,\cdots,X_{N-1}=a_{N-1},X_N=1\right)f\left(a_1,\cdots,a_{N-1},1\right)
\\&=\sum\limits_{a_1,\cdots,a_{N-1}}\p\left(X_1=a_1,\cdots,X_{N-1}=a_{N-1}\right)f\left(a_1,\cdots,a_{N-1},0\right)
\\&+\sum\limits_{a_1,\cdots,a_{N-1}}\p\left(X_1=a_1,\cdots,X_{N-1}=a_{N-1},X_N=1\right)[f\left(a_1,\cdots,a_{N-1},1\right)-f\left(a_1,\cdots,a_{N-1},0\right)]
\\&\le\sum\limits_{a_1,\cdots,a_{N-1}}\p\left(X_1=a_1,\cdots,X_{N-1}=a_{N-1}\right)f\left(a_1,\cdots,a_{N-1},0\right)
\\&+\sum\limits_{a_1,\cdots,a_{N-1}}\p\left(X_1=a_1,\cdots,X_{N-1}=a_{N-1}\right)p[f\left(a_1,\cdots,a_{N-1},1\right)-f\left(a_1,\cdots,a_{N-1},0\right)]
\\&=\left(1-p\right)\sum\limits_{a_1,\cdots,a_{N-1}}\p\left(X_1=a_1,\cdots,X_{N-1}=a_{N-1}\right)f\left(a_1,\cdots,a_{N-1},0\right)
\\&+p\sum\limits_{a_1,\cdots,a_{N-1}}\p\left(X_1=a_1,\cdots,X_{N-1}=a_{N-1}\right)f\left(a_1,\cdots,a_{N-1},1\right)
\\&\triangleq\left(1-p\right)\e f_0\left(X_1,\cdots,X_{N-1}\right)+p\e f_1\left(X_1,\cdots,X_{N-1}\right).
\eae

Since $f_0$ and $f_1$ are both increasing functions on $\{0,1\}^{N-1}$, by the inductive hypothesis we have
\bae
 &\left(1-p\right)\e f_0\left(X_1,\cdots,X_{N-1}\right)+p\e f_1\left(X_1,\cdots,X_{N-1}\right)
 \\&\le\left(1-p\right)\e f_0\left(Y_1,\cdots,Y_{N-1}\right)+p\e f_1\left(Y_1,\cdots,Y_{N-1}\right)
 \\&=\e f\left(Y_1,\cdots,Y_N\right).
 \eae
 Thus we get the result when $n=N$ and the proof is complete.
 \end{proof}

\bibliographystyle{plain}
\bibliography{cite}

\end{document}